# DUAL FORMULATION OF THE UTILITY MAXIMIZATION PROBLEM: THE CASE OF NONSMOOTH UTILITY

By B. Bouchard, N. Touzi and A. Zeghal

*CREST, CREST and CEREMADE*

We study the dual formulation of the utility maximization problem in incomplete markets when the utility function is finitely valued on the whole real line. We extend the existing results in this literature in two directions. First, we allow for nonsmooth utility functions, so as to include the shortfall minimization problems in our framework. Second, we allow for the presence of some given liability or a random endowment. In particular, these results provide a dual formulation of the utility indifference valuation rule.

**1. Introduction.** Given a concave nondecreasing function $U$, finitely valued on the whole real line, we study the dual formulation of the utility maximization problem

$$\sup_{\theta \in \mathcal{H}} EU(X_T^{x,\theta} - B).$$

Here, $X^{x,\theta}$ is the wealth process produced by an initial capital $x$ together with an admissible trading strategy $\theta \in \mathcal{H}$ and $B$ is a given bounded contingent claim, which can also be interpreted as a random endowment. We refer to [17] for an intuitive presentation of the dual problem, although this overview does not address the existence issue.

This problem has been addressed [7] in the context of exponential utility functions. The case of arbitrary smooth utility functions, satisfying the Inada conditions, was studied [18] when $B = 0$. The case of a bounded $B$ was addressed [1] in the presence of transaction costs.

In this article, we focus on the case where the utility function is not assumed to be smooth. Such situations arise naturally in financial markets with transaction costs as argued in [5]. They also appear in many problems in frictionless incomplete markets, such as the shortfall minimization problems studied in [2], [3], [8] and [15], among others.

Our main contribution is the extension of the duality result in [18] and [14] to the above context. In particular, it provides a dual formulation for the Hodges and Neuberger utility-based price; see [9], [1] and [14], among others.

This result is obtained by approximating the utility function by a sequence of utility functions with bounded negative domain. As a by-product, we prove an extension, to the nonsmooth case, of the duality result of [12], which was formulated for utility functions with positive effective domain and $B = 0$. We finally discuss the important issue of the choice of the set of admissible strategies, as addressed in [7] and [19]. We show that the conclusions in [19] extend immediately to our context.

The article is organized as follows. The precise formulation of the problem is presented in Section 2. The main duality results are reported in Section 3 and the discussion on the set of admissible strategies is contained in Section 4. The proofs are collected in the remaining sections.

## 2. Problem formulation.

2.1. *The financial market.* Let $T$ be a finite time horizon and let $(\Omega, \mathcal{F}, P)$ be a complete probability space endowed with a filtration $\mathbb{F} = \{\mathcal{F}_t, 0 \leq t \leq T\}$ satisfying the usual conditions.

The financial market consists of one bank account, with constant price $S^0$, normalized to unity, and $d$ risky assets $S^1, \ldots, S^d$. As usual, there is no loss of generality in normalizing the nonrisky asset price process, since we may always choose it as numeraire under very mild conditions. We denote $S := (S^1, \ldots, S^d)$ the price process of the $d$ risky assets. The vector process $S = \{S_t, 0 \leq t \leq T\}$ is assumed to be a $(0, \infty)^d$-valued semimartingale on the filtered probability space $(\Omega, \mathcal{F}, \mathbb{F}, P)$. Moreover, we assume that the condition

$$(2.1) \qquad \mathcal{M}^e(S) := \{Q \sim P : S \text{ is a } Q\text{-local martingale}\} \neq \varnothing$$

holds. This condition is intimately related to the absence of arbitrage opportunities on the security market; see [6].

A trading strategy $\theta$ is an element of $L(S)$, the set of all $\mathbb{R}^d$-valued predictable processes which are integrable with respect to $S$. In economic words, each component $\theta_t^i$ represents the number of shares of the $i$th risky asset held at time $t$.

Given a trading strategy $\theta \in L(S)$ and initial capital $x \in \mathbb{R}$, it follows from the self-financing condition that the wealth process is defined by

$$X_t^{x,\theta} := x + \int_0^t \theta_r \, dS_r.$$

The possible terminal values of such wealth processes are collected in the set

$$\mathcal{X}(x) := \{X \in L^0 : X = X_T^{x,\theta} \text{ for some } \theta \in L(S)\}.$$

To exclude arbitrage opportunities, it is well known that we need to impose some lower bound on the wealth process. We therefore introduce the subset of $\mathcal{X}(x)$,

$$\mathcal{X}_{\mathrm{b}}(x) := \{X \in \mathcal{X}(x) : \|X^-\|_\infty < \infty\}.$$

2.2. *The utility maximization problem.* Let $U$ be a nonconstant, nondecreasing, concave function defined and finite on the whole real line:

$$\mathrm{dom}(U) := \{x \in \mathbb{R} : |U(x)| < \infty\} = \mathbb{R}.$$

Observe that $U$ is not assumed to be smooth.

In this article, we focus on the problem of maximizing the expected utility from terminal wealth for an agent subject to some liability $B \in L^\infty$. We refer to [14] and [10] for possible extension in the unbounded case. Since existence may fail to hold in $\mathcal{X}_b(x)$ (even in the smooth utility case with $B = 0$), we follow [18] by defining the set $\mathcal{X}_U(x)$ of random variables $X \in L^0$ such that there exists a sequence $X_n \in \mathcal{X}_b(x)$ that satisfies

$$U(X_n - B) \to U(X - B) \qquad \text{in } L^1.$$

We then define the utility maximization problem

$$V(x) := \sup_{X \in \mathcal{X}_U(x)} EU(X - B).$$

Observe that, with this definition, $V(x)$ is also the supremum of the expected terminal wealth over $\mathcal{X}_b(x)$. We conclude this section with some examples of interest in the literature which fit in our framework.

EXAMPLE 2.1 (Smooth utility functions, no liability). When $U$ is continuously differentiable, strictly concave and $B = 0$, the above problem has been addressed in [18]. The particular exponential utility case $U(x) = -e^{-\eta x}$ has been extensively studied in [7] and [11].

EXAMPLE 2.2 (Smooth utility functions with liability). When $U$ is continuously differentiable and strictly concave, the extension to $B \neq 0$ has been performed in [1] and [14]. The main result of this article improves the results of [14] by allowing for a nonsmooth utility function $U$.

EXAMPLE 2.3 (Shortfall utility). Let $\ell$ be a convex nondecreasing function defined on the nonnegative real line. The shortfall minimization problem is defined by

$$\inf_{X \in \mathcal{X}_b(x)} E\ell([B - X]^+).$$

We refer to [2], [3], [8] and [15], among others. Defining $U(x) = -\ell(x^-)$, we see that this problem fits in our framework under mild conditions on $\ell$; see Example 2.4.

2.3. *The dual problem.* Let $\tilde{U}$ be the Legendre–Fenchel transform defined by

$$\tilde{U}(y) := \sup_{x \in \mathbb{R}} U(x) - xy$$

and observe that $\mathrm{dom}(\tilde{U}) \cap (-\infty, 0) = \varnothing$. We assume that the utility function $U$ satisfies

(2.2) $\quad \inf \bigcup_{x \in \mathbb{R}} \partial U(x) = 0 \quad \text{and} \quad r := \sup \bigcup_{x \in \mathbb{R}} \partial U(x) \notin \bigcup_{x \in \mathbb{R}} \partial U(x),$

which can be stated equivalently on $\tilde{U}$ as

(2.3) $\qquad\qquad \mathrm{int}[\mathrm{dom}(\tilde{U})] = (0, r) \quad \text{and} \quad r \notin \mathrm{dom}(\tilde{U}).$

Without loss of generality, we may assume that

(2.4) $\qquad\qquad\qquad\qquad U(0) > 0$

so that $\tilde{U} > 0$. Observe that $r > 0$ since $U$ is not constant. The natural set of dual variables is

(2.5) $\mathbf{Y}_+ := \{(y, Y) \in \mathbb{R}_+ \times L^0_+ : EXY \le xy \text{ for all } x \in \mathbb{R}_+ \text{ and } X \in \mathcal{X}_+(x)\},$

the positive polar of the set of nonnegative elements of $\mathcal{X}(x)$:

(2.6) $\qquad\qquad\qquad \mathcal{X}_+(x) := \mathcal{X}(x) \cap L^0_+.$

However, since we are dealing with a utility function finitely defined on the real line, it turns out that

(2.7) $\qquad\qquad\qquad \tilde{\mathbf{Y}}_+ := \{(y, Y) \in \mathbf{Y}_+ : EY = y\}$

is the appropriate set of dual variables, as was observed in [18]. This set is clearly nonempty because it contains all pairs $(1, Y)$, where $Y = dQ/dP$ with $Q \in \mathcal{M}^e(S)$.

We define the dual problem

$$W(x) := \inf_{(y,Y) \in \tilde{\mathbf{Y}}_+} E[\tilde{U}(Y) + xy - YB].$$

Clearly, we have

(2.8) $$W(x) \geq V(x) \qquad \text{for all } x \in \mathbb{R}.$$

The purpose of this article is to find conditions under which equality holds in the above inequality and to relate the solutions of both problems by the classical Fenchel duality results.

EXAMPLE 2.4 (Back to shortfall utility). In the case of the shortfall utility function $U(x) := -\ell(x^-)$, we directly compute that

$$\tilde{U}(y) := -\inf_{x \geq 0}(\ell(x) - xy).$$

Observe that $\inf \bigcup_{x \in \mathbb{R}} \partial U(x) = 0$ so that this example fits in our framework as long as $\ell$ is not linear near $+\infty$. For instance, for $\ell(x) = x^2$, we compute directly that $\tilde{U}(y) = y^2/4$ and $\text{dom}(\tilde{U}) = [0, \infty)$. However, the case $U(x) = -x^-$ studied in [2] is not covered here.

2.4. *Asymptotic elasticity in the nonsmooth case.* As in [12] and [18], we need conditions on the asymptotic elasticity of the utility function to prove the required duality relationship. In the nonsmooth case, it is argued [5] that these conditions have to be written on the conjugate function $\tilde{U}$. We then define

$$\text{AE}_0(\tilde{U}) := \limsup_{y \searrow 0} \sup_{q \in \partial \tilde{U}(y)} \frac{|q|y}{\tilde{U}(y)} \quad \text{and} \quad \text{AE}_r(\tilde{U}) := \limsup_{y \nearrow r} \sup_{q \in \partial \tilde{U}(y)} \frac{|q|y}{\tilde{U}(y)},$$

where $r$ is the right boundary of the domain of $\tilde{U}$; see (2.3). We show in Lemma 2.2 that the asymptotic elasticity condition $\text{AE}_r(\tilde{U}) < \infty$ together with (2.3) implies that the domain of $\tilde{U}$ is unbounded. We start with the following lemma.

LEMMA 2.1. *Let $f$ be a convex function with $\text{int}[\text{dom}(f)] \subset (0, r)$ for some $r \in (0, \infty] \setminus \text{dom}(f)$. For $k = (k_1, k_2) \in \mathbb{R}_+ \times \mathbb{R}_+$ define*

$$f^k(y) := f(y) - k_1 y + k_2, \qquad y \in \text{dom}(f).$$

*Then:*

(i) *If $f(0+) > 0$, then $\text{AE}_0(f) < \infty \Rightarrow \text{AE}_0(f^k) < \infty$.*

(ii) *If $f(r-) > 0$ and $\liminf_{y \nearrow r} \min \partial \tilde{U}(y) = \infty$, then $\text{AE}_r(f) < \infty \Rightarrow \text{AE}_r(f^k) < \infty$.*

PROOF. (i) Assume that $\text{AE}_0(f) < \infty$. Then there exists a constant $C > 0$ such that, for all sufficiently small $y > 0$ and all $q \in \partial f(y)$, $|q|y \leq Cf(y)$. It follows that, for small $y > 0$,

$$|q - k_1|y \leq |q|y + k_1 y \leq C(f(y) + k_2 - k_1 y) + (C+1)k_1 y \leq C(1 + f^k(y)).$$

Since $f(0+) > 0$, there exists some $\varepsilon > 0$ such that $f^k(y) = f(y) + k_2 - k_1 y > \varepsilon$ for small $y > 0$ and therefore

$$\frac{|q - k_1|y}{f^k(y)} \le C\left(\frac{1}{f^k(y)} + 1\right) \le C\left(\frac{1}{\varepsilon} + 1\right).$$

The result follows.

(ii) Assume that $\mathrm{AE}_r(f) < \infty$. Then there exists a constant $C > 0$ such that, for all $y$ in a neighborhood of $r$ and $q \in \partial f(y)$, $|q|y \le Cf(y)$. This implies

$$qy - k_1 y \le C(f(y) + k_2 - k_1 y) + (C - 1)k_1 y = Cf^k(y) + (C - 1)k_1 y.$$

Since $y > 0$, it follows that

$$q - k_1 \le Cf^k(y)/y + (C - 1)k_1.$$

Since $\liminf_{y \nearrow r} \min \partial f(y) = \infty$ and $q \in \partial f(y)$, we see that, on a neighborhood of $r$, $q - k_1 > 0$, $f^k(y) > 0$ and $f^k(y)/y > \varepsilon$ for some $\varepsilon > 0$. It follows that

$$\frac{|q - k_1|y}{f^k(y)} \le C + \frac{(C-1)k_1}{f^k(y)/y} \le C + (C-1)\frac{k_1}{\varepsilon},$$

which concludes the proof.   □

REMARK 2.1.  Let $U$ be a concave function on $\mathbb{R}$ satisfying (2.2) and let $\tilde{U}$ be the associated Fenchel transform. Then, writing that $-x \in \partial \tilde{U}(y) \Rightarrow y \in \partial U(x)$ (see, e.g., [16]) implies that $\liminf_{y \nearrow r} \min \partial \tilde{U}(y) = \infty$. In view of (2.3), we see that Lemma 2.1 applies for $f = \tilde{U}$. For later purposes, observe that this implies that $\tilde{U}$ is nondecreasing near $r \in (0, \infty]$.

LEMMA 2.2.  *Assume that the conjugate function satisfies* (2.3) *as well as the asymptotic elasticity condition* $\mathrm{AE}_r(\tilde{U}) < \infty$. *Then* $r = +\infty$.

PROOF.  We assume that $r < \infty$ and work toward a contradiction.

STEP 1.  We first prove that we can assume w.l.o.g. that $\tilde{U}$ is positive and nondecreasing near $r$. To see this, define $U^k(x) = U(x - k_1) + k_2$ for $k = (k_1, k_2) \in \mathbb{R}_+ \times \mathbb{R}_+$. From (2.3), observe that we can choose $k$ such that $U^k(0) > 0$ and $\max \partial U^k(0) < r$, so that $\tilde{U}^k$ is positive and nondecreasing near $r$. Using Lemma 2.1 and Remark 2.1, we can then reduce the statement of the lemma to $\tilde{U}^k(y) = \tilde{U}(y) - k_1 y + k_2$ since $\mathrm{dom}(\tilde{U}) = \mathrm{dom}(\tilde{U}^k)$ and $\mathrm{AE}_r(\tilde{U}) < \infty$ implies $\mathrm{AE}_r(\tilde{U}^k) < \infty$.

STEP 2. From Step 1, we can assume that $\tilde{U}$ is positive and nondecreasing near $r$. Now observe that $\mathrm{AE}_r(\tilde{U}) < \infty$ implies the existence of some constant $C$ such that $\max \partial \tilde{U}(y)/\tilde{U}(y) \leq C$ for all $y \in [r', r)$ for some $r' < r$. Then, for all $y \in [r', r)$, $\tilde{U}(y) \leq \alpha e^{Cy}$ for some real $\alpha$. Since $r < \infty$, this implies that $\tilde{U}(r-) < \infty$. We conclude the proof by observing that any $x' \in \partial \tilde{U}(r-)$ satisfies $r \in \partial U(-x')$ by the classical connection between the gradients of $U$ and $\tilde{U}$; see, for example, [16]. This contradicts (2.3). □

In view of this result, we rewrite (2.3) as

$$(2.9) \qquad \mathrm{int}[\mathrm{dom}(\tilde{U})] = (0, \infty).$$

The following result is an extension to the nonsmooth case of the implications of the asymptotic elasticity conditions derived in [18]. We postpone its proof to Section 7.

LEMMA 2.3. *Let $f$ be a positive convex function, with $\mathrm{cl}[\mathrm{dom}(f)] = \mathbb{R}_+$. Assume further that $f$ is nonincreasing near $0$, nondecreasing near $\infty$ and satisfies the asymptotic elasticity conditions*

$$(2.10) \qquad \mathrm{AE}_0(f) < \infty \quad and \quad \mathrm{AE}_\infty(f) < \infty.$$

*Then for all $0 < \mu_0 < \mu_1 < \infty$, there exists a constant $C > 0$ such that:*

   (i) *$f(\mu y) \leq C f(y)$ for all $\mu \in [\mu_0, \mu_1]$ and $y > 0$;*
   (ii) *$y|q| \leq C f(y)$ for all $y > 0$ and $q \in \partial f(y)$.*

## 3. The main result.

### 3.1. *Utility functions with unbounded domain.*

REMARK 3.1. Up to now, we have not assumed that $S$ is locally bounded. In turns out that this technical assumption is not needed for our result. However, as pointed out in Remark 2.6 of [18], the set of strategies $\mathcal{X}_U$ may not be adapted when $S$ is not locally bounded. More precisely, we can construct easy examples where the primal problem has a natural solution outside $\mathcal{X}_U$ and the restriction of the strategies to $\mathcal{X}_U$ leads to a zero investment strategy as an optimal solution, which makes no sense from an economic point of view. For instance, set $B = 0$ and consider a market with one risky asset $S^1$ such that $S^1 = 1$ on $[0, T)$ and $S_T^1$ is normally distributed (assuming now that prices can be negative), that is, $S^1$ jumps at $T$. Then, it is easily checked that $\mathcal{X}_b(x) = \{x\}$ and therefore $V(x) = U(x)$, that is, the optimal strategy in $\mathcal{X}_U(x)$ is $X_* = x$. Assuming that $U$ is strictly concave and smooth. Since $\mathcal{X}_+(r) = \{r\}$ for $r \geq 0$, we see that $(y_*, Y_*) = (U'(x), U'(x)) \in \tilde{\mathbf{Y}}_+$. Since $\tilde{U}(U'(x)) + x U'(x) = U(x)$, we also see that the usual duality holds and that $(y_*, Y_*)$ is optimal for $W(x)$, and we easily check that all the requirements of Theorem 3.1 are satisfied, except that $Y_*/y_* = 1$ does not define a local martingale measure if $ES_T^1 \neq 1$.

In view of this remark, we assume in this subsection that $S$ is locally bounded. This will prevent the above described phenomenon.

REMARK 3.2.  Define the sequence of stopping times $\tau_n := \inf\{t \geq 0 : |S_t| > n\}$. Since $S$ is locally bounded, we have $S_{\tau_n}^i \in \mathcal{X}_b(S_0^i)$ and $-S_{\tau_n}^i \in \mathcal{X}_b(-S_0^i)$. By definition of $\tilde{\mathbf{Y}}_+$, we deduce that, for each $(y, Y) \in \tilde{\mathbf{Y}}_+$ with $y > 0$, the measure $Q := (Y/y) \cdot P \in \mathcal{M}^a(S)$, the set of all local martingale measures for $S$ which are absolutely continuous with respect to $P$.

THEOREM 3.1.   *Let $U$ be a nonconstant concave nondecreasing function, finitely valued in $\mathbb{R}$, satisfying* (2.4) *and such that the associated Fenchel transform $\tilde{U}$ satisfies* (2.9) *as well as the asymptotic elasticity conditions* (2.10). *Given some bounded contingent claim $B$, consider the optimization problems*

$$V(x) := \sup_{X \in \mathcal{X}_U(x)} EU(X - B) \quad and \quad W(x) := \inf_{(y,Y) \in \tilde{\mathbf{Y}}_+} E[\tilde{U}(Y) + yx - YB].$$

*Assume further that $W(x) < \infty$ for some $x \in \mathbb{R}$. Then:*

(i) *Existence holds for the dual problem $W(x)$, that is,*

$$W(x) = E[\tilde{U}(Y_*) - Y_*B + xy_*] \qquad for\ some\ (y_*, Y_*) \in \tilde{\mathbf{Y}}_+.$$

*Moreover, if $y_* > 0$, then $Q_* = \frac{Y_*}{y_*} \cdot P \in \mathcal{M}^a(S)$.*

(ii) *Existence holds for the portfolio optimization problem $V(x)$, that is,*

$$V(x) = E[U(X_* - B)] \qquad for\ some\ X_* \in \mathcal{X}_U(x).$$

(iii) *The above solutions are related by*

$$X_* \in B - \partial \tilde{U}(Y_*), \qquad P\text{-}a.s. \quad and \quad E[X_*Y_*] = xy_*,$$

*so that the duality relationship $V(x) = W(x)$ holds.*

(iv) *If $Y_* > 0$, $P$-a.s., then $X_* = X_T^{x,\theta}$ for some $\theta \in L(S)$, where $X^{x,\theta}$ is a uniformly integrable martingale under the measure $Q_* := \frac{Y_*}{y_*} \cdot P \in \mathcal{M}^e(S)$.*

The proof of this result is reported in Section 5.

REMARK 3.3.   It is immediately checked that

$$W(x) < \infty \text{ for some } x \in \mathbb{R} \begin{cases} \text{if and only if } E[\tilde{U}(Y)] < \infty \text{ for some } (y,Y) \in \tilde{\mathbf{Y}}_+, \\ \text{if and only if } W(x) < \infty \text{ for all } x \in \mathbb{R}. \end{cases}$$

We next focus on the attainability issue of Theorem 3.1(iv). Clearly, since $\tilde{U}(0) = U(\infty)$, it follows from Remark 3.3 that $Y_* > 0$ whenever $U(\infty) = \infty$. More generally, we prove the following sufficient condition in Section 5.

PROPOSITION 3.1.  *Assume that $L := \inf\{l : U(l) = U(\infty)\} = \infty$. In the context of Theorem 3.1, assume further that $\tilde{\mathbf{Y}}_+$ contains some $(\bar{y}, \bar{Y})$ satisfying*

$$E\tilde{U}(\bar{Y}) < \infty \quad and \quad \bar{Y} > 0, \qquad \text{$P$-a.s.}$$

*Then $Y_* > 0$, $P$-a.s.*

REMARK 3.4.  We now discuss the uniqueness issue when the utility function $U$ is strictly concave. Observe that $\mathcal{X}_U$ is a priori not convex. However, we see in this remark that this property holds if we restrict to the set of optimal strategies, thus providing uniqueness. Let $X_*^1$ and $X_*^2$ be two solutions of the utility maximization problem and let $X_n^1, X_n^2 \in \mathcal{X}_b(x)$ be such that $U(X_n^i - B) \to U(X_*^i - B)$ in $L^1$, $i \in \{1, 2\}$. Since $U$ is increasing, we see that, possibly after passing to subsequences, $X_n^i \to X_*^i$, $P$-a.s., $i \in \{1, 2\}$. Since, for all $\lambda \in (0, 1)$,

$$U(\lambda X_n^1 + (1-\lambda)X_n^2 - B) \geq \lambda U(X_n^1 - B) + (1-\lambda)U(X_n^2 - B)$$
$$\to \lambda U(X_*^1 - B) + (1-\lambda)U(X_*^2 - B) \qquad \text{in } L^1$$

and $\lambda X_n^1 + (1-\lambda)X_n^2 \in \mathcal{X}_b(x)$, it follows that

$$V(x) = \lim_{n \to \infty} EU(\lambda X_n^1 + (1-\lambda)X_n^2 - B)$$

and

$$U(\lambda X_n^1 + (1-\lambda)X_n^2 - B) \to U(\lambda X_*^1 + (1-\lambda)X_*^2 - B) \qquad \text{in } L^1,$$

$$V(x) = EU(\lambda X_*^1 + (1-\lambda)X_*^2 - B) = \lambda EU(X_*^1 - B) + (1-\lambda)EU(X_*^2 - B).$$

It follows that, in the case where $U$ is strictly concave, there is a unique solution to the utility maximization problem. However, if $U$ is not smooth, the Fenchel transform $\tilde{U}$ is not strictly convex and uniqueness in the dual problem is not guaranteed. We continue this discussion in Remark 4.2. We thank an anonymous referee for pointing out this important issue.

To prove Theorem 3.1, we use the methodology of [18], which consists of approximating $U$ by utility functions $U_n$ that have a domain bounded from below. Set

$$U_n := U \qquad \text{on } \mathrm{dom}(U_n) := (-n, \infty) \qquad \text{for } n \geq 2\|B\|_\infty,$$

so that $U_n$ converges to $U$ and $\mathrm{dom}(U_n)$ is bounded from below. Let $\tilde{U}_n$ be the associated Fenchel transform

$$\tilde{U}_n(y) := \sup_{x \in \mathbb{R}}(U_n(x) - yx).$$

Observe that our approximating utility functions are nonsmooth and that

$$(3.1) \quad U_n = U \qquad \text{on } \mathrm{dom}(U_n) \quad \text{and} \quad \tilde{U}_n = \tilde{U} \qquad \text{on } \partial U_n \, \mathrm{dom}(U_n).$$

We follow [1] by defining

$$x_n := x + \frac{n}{2} \quad \text{and} \quad B_n := B + \frac{n}{2},$$

together with the corresponding approximating optimization problems

$$V_n(x) := \sup_{X \in \mathcal{C}(x_n)} EU_n(X - B_n),$$

$$W_n(x) := \inf_{(y,Y) \in \mathbf{Y}_+} E\tilde{U}_n(Y) - YB_n + x_n y,$$

where $\mathbf{Y}_+$ is defined in (2.5) and

$$\mathcal{C}(x) := \{X \in L^0_+ - L^\infty : EXY \le xy \text{ for all } (y,Y) \in \mathbf{Y}_+\}.$$

The reason for introducing the sequences $(x_n)_n$ and $(B_n)_n$ appears in Lemma 5.4.

REMARK 3.5.  Since $\mathbf{Y}_+$ contains all pairs $(1, dQ/dP)$ for $Q \in \mathcal{M}^e(S)$, it follows from the classical dual formulation of the superreplication problem that

$$\mathcal{C}(x) \subset \{X \in L^0_+ - L^\infty : X \le X^s \text{ } P\text{-a.s. for some } X^s \in \mathcal{X}_\mathrm{b}(x)\};$$

that is, all contingent claims in $\mathcal{C}(x)$ can be superreplicated starting from the initial capital $x$. By definition of $\mathbf{Y}_+$, the reverse inclusion holds for nonnegative contingent claims, so that

$$\mathcal{C}(x) \cap L^0_+ = \{X \in L^0_+ : X \le X^s \text{ for some } X^s \in \mathcal{X}_+(x)\}.$$

The first step in the proof of Theorem 3.1 is to establish existence for the above approximating control problems as well as the duality connection between them. This is the main object of the following subsection.

3.2. *Utility functions with bounded negative domain.*  We now concentrate on the case where the utility function has a domain which is bounded from below.

THEOREM 3.2.   *Let $\beta \geq 0$ be an arbitrary constant and consider some contingent claim $B$ with $\|B\|_\infty \leq \beta$. Let $U$ be a nonconstant concave non-decreasing function with*

$$\text{cl}[\text{dom}(U)] = [-2\beta, \infty), \qquad U(\infty) > 0, \qquad \text{cl}[\text{dom}(\tilde{U})] = \mathbb{R}_+$$

*and satisfying the asymptotic elasticity condition $\text{AE}_0(\tilde{U}) < \infty$. Consider the optimization problems*

$$V(x) := \sup_{X \in \mathcal{C}(x)} EU(X - B) \quad and \quad W(x) := \inf_{(y,Y) \in \mathbf{Y}_+} E[\tilde{U}(Y) + yx - YB].$$

*Assume that $W(x) < \infty$ for some $x > 0$. Then:*

(i)  *Existence holds for the dual problem $W(x)$, that is,*

$$W(x) = E[\tilde{U}(Y_*) + y_* x - Y_* B] \qquad for \; some \; (y_*, Y_*) \in \mathbf{Y}_+.$$

(ii)  *Existence holds for the optimization problem $V(x)$, that is,*

$$V(x) = E[U(X_* - B)] \qquad for \; some \; X_* \in \mathcal{C}(x) \; such \; that \; X_* - B \geq -2\beta.$$

*Moreover, if $X_* \geq 0$, then $X_* \in \mathcal{X}_+(x)$.*

(iii)  *The above solutions are related by*

$$X_* \in B - \partial\tilde{U}(Y_*), \qquad P\text{-a.s.} \quad and \quad E[X_* Y_*] = xy_*,$$

*so that the duality relationship $V(x) = W(x)$ holds.*

The proof is postponed to Section 6.

REMARK 3.6.   The technical assumption $U(\infty) > 0$ can clearly be relaxed by adding a constant to $U$.

REMARK 3.7.   Corollary 6.3 states that the solution of the dual problem, introduced in Theorem 3.2, satisfies $Y_* > 0$ $P$-a.s. whenever $L := \inf\{l : U(l) = U(\infty)\} = \infty$.

REMARK 3.8.   As in Remark 3.4, we assume that $U$ is strictly concave, so that the solution to the utility maximization problem is unique. Recalling that, for all $(X, y, Y) \in \mathcal{C}(x) \times \mathbf{Y}_+$, $EXY \leq xy$, we see by similar arguments as in Remark 4.2 that uniqueness holds in the dual problem outside of the set where $\tilde{U}'$ is constant.

REMARK 3.9.   Let us specialize the discussion of Theorem 3.2 to the case $B = 0$.

1. First let $\beta = 0$. Then, obviously, $X_*$ is nonnegative and therefore

$$V(x) = \sup_{X \in \mathcal{X}_+(x)} EU(X) = EU(X_*).$$

We are in the context of the portfolio optimization problem of [12], except that the utility function is not assumed to be smooth. Hence, Theorem 3.2 extends the corresponding results to the nonsmooth utility case. It is also easy to check that we have the additional result

$$W(x) = \inf_{y>0} \inf_{Q \in \mathcal{M}^e(S)} E\Big[\tilde{U}\Big(y\frac{dQ}{dP}\Big) + yx\Big]$$

by the same arguments as in [12].

2. For $\beta > 0$ and $x > -2\beta$, the same argument as in [18], Section 2, shows that existence holds for the problem

$$\sup_{X \in \mathcal{X}_b(x)} EU(X)$$

and that the solution $X_*$ of the above problem is related to the solution $\bar{X}_*$ of the problem defined on the utility function $U(\cdot - 2\beta)$ with initial wealth $\bar{x} = x + 2\beta$ by $X_* = \bar{X}_* - 2\beta$.

3. Because of the connection between $\|B\|_\infty$ and the domain of $U$, and the nature of the set of primal variables $\mathcal{C}(x)$, Theorem 3.2 does not compare to [4] and [10].

**4. Complements on the set of admissible strategies in the unbounded domain case.** Following [18], we now consider alternative sets of admissible strategies for the problem of Section 3.1. In view of Remark 3.1, we assume that $S$ is locally bounded. Recall from Remark 3.2 that, under this condition, for each $(y, Y) \in \tilde{\mathbf{Y}}_+$ with $y > 0$, the measure $Q := (Y/y) \cdot P \in \mathcal{M}^a(S)$.

Let $x \in \mathbb{R}$ be some fixed initial capital and assume that the conditions of Theorem 3.1 hold, so that solutions $X_*$ of $V(x)$ and $(y_*, Y_*)$ of $W(x)$ do exist and satisfy the conditions of the theorem. Then, if $y_* > 0$, the induced measure

$$Q_* := \frac{Y_*}{y_*} \cdot P \in \mathcal{M}^a_{\tilde{U}}(S) := \Big\{ Q \in \mathcal{M}^a(S) : E\tilde{U}\Big(\frac{dQ}{dP}\Big) < \infty \Big\}.$$

Throughout this section, we assume that $Y_*$ satisfies the additional condition

$$Y_* > 0, \qquad P\text{-a.s.},$$

so that $y_* > 0$, $Q_* \in \mathcal{M}^e(S)$ and

$$W(x) = \inf_{\substack{y>0 \\ Q \in \mathcal{M}^e(S)}} E\Big[\tilde{U}\Big(y\frac{dQ}{dP}\Big) - By\frac{dQ}{dP}\Big] + xy.$$

The measure $Q_*$ is the so-called minimal local martingale measure associated to the problem $\tilde{V}(y_*)$, where

$$\tilde{V}(y) := \inf_{(y,Y) \in \tilde{\mathbf{Y}}_+} E[\tilde{U}(Y) - YB].$$

Under the assumption $Y_* > 0$, we also know from Theorem 3.1 that $X_* = X_T^{x,\theta^*}$ for some $\theta^* \in L(S)$.

A simple restatement of Theorem 3.1(iii) and (iv) reveals that the wealth process $X^{x,\theta_*}$ is a uniformly integrable martingale under $Q_*$,

and

$$V(x) = \inf_{y>0} \tilde{V}(y) + xy \qquad \text{so that } x \in -\partial\tilde{V}(y_*),$$

where we used the (obvious) convexity of $\tilde{V}$. The following sets of strategies were studied in [7] and [19]:

$\mathcal{H}_1(x) := \{\theta \in L(S) : U(X_T^{x,\theta} - B) \in L^1 \text{ and } X^{x,\theta} \text{ is a } Q_*\text{-supermartingale}\}$,

$\mathcal{H}_1'(x) := \{\theta \in \mathcal{H}_1(x) : X^{x,\theta} \text{ is a } Q_*\text{-martingale}\}$,

$\mathcal{H}_2(x) := \{\theta \in \mathcal{H}_1(x) : X^{x,\theta} \text{ is a supermartingale under all } Q \in \mathcal{M}_{\tilde{U}}^a(S)\}$.

We now have the following extension of [19] to the nonsmooth utility context of this article.

1. Clearly, since $EU(X_T^{x,\theta} - B) \leq W(x)$ for all $\theta \in \mathcal{H}_1(x)$ and $\theta_* \in \mathcal{H}_1'(x) \subset \mathcal{H}_1(x)$, it follows that

$$V(x) = \max_{\theta \in \mathcal{H}_1(x)} EU(X_T^{x,\theta} - B) = \max_{\theta \in \mathcal{H}_1'(x)} EU(X_T^{x,\theta} - B).$$

2. Also, observe that $\mathcal{X}_b(x) \subset \{X_T^{x,\theta} : \theta \in \mathcal{H}_2(x)\}$. Therefore,

$$V(x) \leq \sup_{\theta \in \mathcal{H}_2(x)} EU(X_T^{x,\theta} - B)$$

$$\leq \inf_{y>0, Q \in \mathcal{M}_{\tilde{U}}^a(S)} E\left[\tilde{U}\left(y\frac{dQ}{dP}\right) - y\frac{dQ}{dP}B + yx\right]$$

$$\leq E\left[\tilde{U}\left(y_*\frac{dQ_*}{dP}\right) - y_*\frac{dQ_*}{dP}B + y_*x\right] = W(x) = V(x).$$

Hence equality holds in all the above inequalities. In particular, this proves that

$$V(x) = \sup_{\theta \in \mathcal{H}_2(x)} EU(X_T^{x,\theta} - B).$$

3. We now prove that $\theta_* \in \mathcal{H}_2(x)$ so that

$$(4.1) \qquad V(x) = \max_{\theta \in \mathcal{H}_2(x)} EU(X_T^{x,\theta} - B).$$

Let $F$ be the conjugate of the function $x \mapsto U(x - \|B\|_\infty)$, that is,

$$F : y \mapsto \tilde{U}(y) - y\|B\|_\infty.$$

Arguing as in Lemma 5.1, we may assume without loss of generality that

$$(4.2) \qquad F(0) > 0, \qquad F \text{ is nonincreasing near } 0,$$

$$(4.3) \qquad \mathrm{AE}_0(F) < \infty \quad \text{and} \quad \mathrm{AE}_\infty(F) < \infty.$$

Notice that, by Remark 2.1 and (2.9), $F$ is clearly nondecreasing near $+\infty$.

To see that (4.1) holds, it suffices to prove that the conjugate function $\tilde{V}$ inherits the asymptotic elasticity conditions $\mathrm{AE}_0(\tilde{V}) < \infty$ and $\mathrm{AE}_{+\infty}(\tilde{V}) < \infty$ from the function $\tilde{U}$. In view of the above assumptions (4.2), we need to show that

$$(4.4) \qquad \begin{aligned} &\text{for all } 0 < \mu_0 < \mu_1, \text{ there exists some } C > 0, \\ &\tilde{V}(\lambda y) \le C\tilde{V}(y) \qquad \text{for all } \lambda \in [\mu_0, \mu_1] \text{ and } y > 0. \end{aligned}$$

With this property of $\tilde{V}$, the proof of Proposition 2.2 in [19] applies immediately to the nonsmooth case.

The characterization of the asymptotic elasticity conditions of Lemma 2.3 holds for $F$ by (4.2), (4.3) and the fact that it is nondecreasing near $+\infty$. Let $(y, Y^\varepsilon) \in \bar{\mathbf{Y}}_+$ be such that

$$E[\tilde{U}(Y^\varepsilon) - Y^\varepsilon B] \le \tilde{V}(y) + \varepsilon.$$

Fix $0 < \mu_0 < \mu_1$. Then, by Lemma 2.3 and (4.3), there exists some $C > \mu_1$ such that for all $\lambda \in [\mu_0, \mu_1]$,

$$\tilde{V}(y) + \varepsilon \ge E[F(Y^\varepsilon) + Y^\varepsilon(\|B\|_\infty - B)] \ge \frac{1}{C} E[F(\lambda Y^\varepsilon) + Y^\varepsilon(\|B\|_\infty - B)]$$

$$= \frac{1}{C} E[\tilde{U}(\lambda Y^\varepsilon) - \lambda Y^\varepsilon B] + \left(1 - \frac{\lambda}{C}\right) E[Y^\varepsilon(\|B\|_\infty - B)]$$

$$\ge \frac{1}{C} E[\tilde{U}(\lambda Y^\varepsilon) - \lambda Y^\varepsilon B] \ge \frac{1}{C} \tilde{V}(\lambda y)$$

and (4.4) follows by arbitrariness of $\varepsilon > 0$.

REMARK 4.1.   It is known from [19] that considering sets of admissible strategies such as

$$\{\theta \in L(S) : X^{x,\theta} \text{ is a } Q\text{-supermartingale (resp. martingale)}$$

$$\text{under some } Q \in \mathcal{M}^e(S)\}$$

may lead to paradoxical results from an economic point of view. They are therefore not discussed in this article.

REMARK 4.2. We continue the discussion on the uniqueness issue of Remark 3.4. It follows from the above analysis that if $S$ is locally bounded and $y_* > 0$, then we are reduced to considering the sets $\mathcal{H}_2(x)$ for the primal problem and $\mathcal{M}_{\tilde{U}}^a(S)$ for the dual problem. Recall from Remark 3.4 that if $U$ is strictly concave, then uniqueness holds in the utility maximization problem. Then, writing $E[(dQ/dP)X^{x,\theta_*}] \leq x$ for all $Q \in \mathcal{M}_{\tilde{U}}^a(S)$, we see that a necessary and sufficient condition for $Q$ to be optimal for the dual problem is that

$$E\left[\frac{dQ}{dP}X_T^{x,\theta_*}\right] = x \quad \text{and} \quad X_T^{x,\theta_*} \in -\partial\tilde{U}\left(y_*\frac{dQ}{dP}\right).$$

It follows that if $U$ is strictly concave and therefore $\tilde{U}$ is continuously differentiable, the optimum for the dual problem is unique outside of the set where $\tilde{U}'$ is constant, that is, $\{y \geq 0 : y \in \partial U(x)$, for some $x$ where $U$ is not differentiable$\}$.

**5. Proofs for the unbounded negative domain case.** In this section, we report the proofs of Theorem 3.1 and Proposition 3.1. We split the proof of Theorem 3.1 into different lemmas. We start by a convenient reduction of the problem.

LEMMA 5.1. *Suppose that statements* (i)–(iv) *of Theorem* 3.1 *hold for* $x > \|B\|_\infty$ *and* $\tilde{U}$ *nonincreasing near* 0. *Then Theorem* 3.1 *holds.*

PROOF. First notice from (2.9) and $U(0) > -\infty$, that for all sufficiently large $k = (k_1, k_2)$, the shifted utility function $U^k : z \in \mathbb{R} \mapsto U(z - k_1) + k_2$ satisfies $\max \partial U^k(0) > 0$ and $U^k(0) > 0$. It follows that the associated Fenchel transform function $\tilde{U}^k$ is positive and, by the classical connection between the gradients $\partial U^k$ and $\partial \tilde{U}^k$ (see, e.g., [16]), that $\tilde{U}^k$ is nonincreasing near 0.

Now, choose $k$ so that the additional condition $x_k := x + k_1 > \|B\|_\infty$ holds. By Lemma 2.1 and Remark 2.1, $\tilde{U}^k$ satisfies the asymptotic elasticity condition of Theorem 3.1; see (5.1). By assumption of the lemma, it follows that Theorem 3.1 holds for the problems

$$V^k(x_k) := \sup_{X \in \mathcal{X}_U(x_k)} EU^kt(X - B)$$

and

$$W^k(x_k) := \inf_{(y,Y)\in\tilde{\mathbf{Y}}_+} E(\tilde{U}^k(Y) - YB) + yx_k.$$

We denote by $(y_*^k, Y_*^k)$ (resp. $X_*^k$) the solution of the problem $W^k(x_k)$ [resp. $V^k(x_k)$]. Observing that for $y \geq 0$,

$$(5.1) \qquad -\partial \tilde{U}^k(y) = -\partial \tilde{U}(y) + k_1, \qquad \tilde{U}^k(y) = \tilde{U}(y) - yk_1 + k_2,$$

it is easily checked that $(y_*, Y_*) := (y_*^k, Y_*^k)$ (resp. $X_* := X_*^k - k_1$) is optimal for the problem $W(x)$ [resp. $V(x)$] and that these quantities satisfy all the statements of Theorem 3.1. □

In view of this result, we assume from now on that

$$x > \|B\|_\infty, \qquad \tilde{U} \text{ is positive and nonincreasing near } 0.$$

We recall from Remark 2.1 and (2.9) that

$$\tilde{U} \text{ is nondecreasing near } +\infty,$$

so that the conditions of Lemma 2.3 hold for $\tilde{U}$.

REMARK 5.1. We isolate the following arguments which will be used repeatedly.

(i) Since $\mathcal{X}_+(1)$ contains the constant random variable 1, we have

$$(5.2) \qquad\qquad EY \leq y \qquad \text{for all } (y, Y) \in \mathbf{Y}_+$$

and, for all constant $M > 0$,

$$\text{the family } \{(y, Y) \in \mathbf{Y}_+ : |y| \leq M\} \text{ is bounded in } L^1(P).$$

(ii) Then, for any sequence $(y_n, Y_n)_n \subset \mathbf{Y}_+$ with bounded $(y_n)_n$, it follows from the Komlòs lemma together with the convexity of $\mathbf{Y}_+$ and Fatou's lemma that

$$\text{there is a sequence } (\tilde{y}_n, \tilde{Y}_n) \in \text{conv}\{(y_k, Y_k), k \geq n\}$$

$$\text{such that } P\text{-a.s.} (\tilde{y}_n, \tilde{Y}_n) \rightarrow (\tilde{y}, \tilde{Y}) \in \mathbf{Y}_+.$$

We now apply Theorem 3.2 to the approximating nonsmooth utility function $U_n$ for some $n \geq 2\|B\|_\infty$. Obviously, $\text{AE}_0(\tilde{U}_n) = \text{AE}_0(\tilde{U}) < \infty$ by (3.1). We need to check only that $W_n(x) < \infty$. In view of Remark 3.3, this is a consequence of the following lemma.

LEMMA 5.2. *The sequence* $(W_n(x))_n$ *is nondecreasing and bounded from above by* $W(x)$.

PROOF. Fix $m > n \in \mathbb{N}$ and consider some $(y, Y) \in \mathbf{Y}_+$. Since $\{\tilde{U}_n\}$ is increasing and $y \geq EY$, we obtain

$$E[\tilde{U}_n(Y) + yx_n - YB_n] \leq E[\tilde{U}_m(Y) + yx_n - YB_n] + \frac{m-n}{2}(y - EY)$$
$$= E[\tilde{U}_m(Y) + yx_m - YB_m].$$

It follows that $(W_n(x))_n$ is nondecreasing. Now fix $(y, Y) \in \tilde{\mathbf{Y}}_+$ and $n \in \mathbb{N}$. Since $\tilde{U}_n \leq \tilde{U}$,

$$E[\tilde{U}_n(Y) + yx_n - YB_n] \leq E[\tilde{U}(Y) + yx - YB] + \frac{n}{2}(y - EY).$$

The required result follows from the fact that $EY = y$ and $\tilde{\mathbf{Y}}_+ \subset \mathbf{Y}_+$. $\square$

We are then in the context of Theorem 3.2. Throughout this section, we denote by $(y_n, Y_n) \in \mathbf{Y}_+$ a solution of problem $W_n(x)$ and by $X_n \in \mathcal{C}(x_n)$ a solution of problem $V_n(x)$ that satisfy the assertions of Theorem 3.2. We recall the connection between these solutions. From (3.1), it follows that

(5.3) $\quad W_n(x) = E[\tilde{U}(Y_n) + x_n y_n - Y_n B_n] = V_n(x) = E[U(X_n - B_n)],$

(5.4) $\quad X_n \in B_n - \partial \tilde{U}_n(Y_n) = B_n - \partial \tilde{U}(Y_n) \quad$ and $\quad E[X_n Y_n] = x_n y_n.$

By Remark 3.5, there exist some $X_n^s \in \mathcal{X}_b(x_n)$ that satisfy $X_n^s \geq X_n$, $P$-a.s. We denote by $V_n^s(x)$ the associated expected utility:

$$V_n^s(x) := EU_n(X_n^s - B_n) = EU(X_n^s - B_n).$$

Observing that $X_n^s - n/2 \in \mathcal{X}_b(x)$, we directly see that

(5.5) $$V_n(x) \leq V_n^s(x) \leq V(x).$$

The following result follows from the same argument as in Step 2 of [18].

LEMMA 5.3. *The sequence $(Y_n)_n$ is uniformly integrable.*

The next result completes the proof of Theorem 3.1(i) and prepares for the proof of the remaining items.

LEMMA 5.4. (i) *There is a sequence $(\hat{y}_n, \hat{Y}_n, \hat{J}_n) \in \mathrm{conv}\{(y_k, Y_k, \tilde{U}(Y_k)), k \geq n\}$ such that*

$$(\hat{y}_n, \hat{Y}_n, \hat{J}_n) \to (y_*, Y_*, \tilde{U}(Y_*)) \in \tilde{\mathbf{Y}}_+ \times L^1(\mathbb{R}_+), \qquad P\text{-a.s. and in } L^1.$$

(ii) *$(y_*, Y_*)$ is optimal for $W(x)$, that is, $(y_*, Y_*) \in \tilde{\mathbf{Y}}_+$, and $E[\tilde{U}(Y_*) + y_* x - Y_* B] = W(x).$*

(iii) $V_n(x) = W_n(x) \uparrow W(x) = V(x) < \infty$ and $V_n^s(x) \to V(x)$.

PROOF.

STEP 1. By (5.2), (5.3), Lemma 5.2 and the positivity of $\tilde{U}$, it follows that

$$\infty > W(x) \geq (x - \|B\|_\infty)y_n + \frac{n}{2}(y_n - EY_n).$$

This proves that $y_n \to y_* \geq 0$ and $y_n - EY_n \to 0$ along some subsequence, as $x - \|B\|_\infty > 0$, $y_n \geq 0$ and $y_n - EY_n \geq 0$. The existence of a sequence $(\hat{y}_n, \hat{Y}_n) \in \mathrm{conv}\{(y_k, Y_k), k \geq n\}$, which converges $P$-a.s. to $(y_*, Y_*) \in \mathbf{Y}_+$, follows from Remark 5.1(ii). From Lemma 5.3, the convergence of $\hat{Y}_n$ to $Y_*$ holds in $L^1$ and therefore $EY_* = y_*$, proving that $(y_*, Y_*) \in \bar{\mathbf{Y}}_+$.

STEP 2. Let $C$ be such that for all $n \geq 2\|B\|_\infty$,

$$\tilde{U}_n(Y_n) - Y_n B \geq U_n(-B) \geq -C > -\infty.$$

Let $(\mu_n^k)_{n,k}$ denote the coefficients of the convex combination that defines the sequence $(\hat{y}_n, \hat{Y}_n)_n$. Using Fatou's lemma, the inequality $y_k \geq EY_k$, Step 1 and (3.1), we get

$$
\begin{aligned}
E(\tilde{U}&(Y_*) + y_* x - Y_* B) \\
&\leq E\left(\liminf_{n \to \infty} \sum_{k \geq n} \mu_n^k(\tilde{U}(Y_k) + y_k x - Y_k B)\right) \\
&\leq \liminf_{n \to \infty} E\left(\sum_{k \geq n} \mu_n^k(\tilde{U}_k(Y_k) + y_k x - Y_k B)\right) \\
&\leq \liminf_{n \to \infty} \sum_{k \geq n} \mu_n^k\left(E[\tilde{U}_k(Y_k) + y_k x - Y_k B] + \frac{k}{2}(y_k - EY_k)\right) \\
&= \liminf_{n \to \infty} \sum_{k \geq n} \mu_n^k W_k(x) \leq W(x) < \infty.
\end{aligned}
$$

Since $(y_*, Y_*) \in \bar{\mathbf{Y}}_+$, it is optimal for $W(x)$. By Lemma 5.2 and (5.3), it follows that

$$(5.7) \quad E(\tilde{U}(Y_*) + y_* x - Y_* B) = W(x) = \lim_{n \to \infty} \uparrow W_n(x) = \lim_{n \to \infty} \uparrow V_n(x).$$

STEP 3. The above argument also proves that $\sup_n E \sum_{k \geq n} \mu_n^k \tilde{U}(Y_k) = \sup_n E \sum_{k \geq n} \mu_n^k |\tilde{U}(Y_k)| < \infty$. We can, therefore, find a sequence $\hat{J}_k \in \mathrm{conv}\{\sum_{k \geq l} \mu_l^k \tilde{U}(Y_k), l \geq n\}$ which converges $P$-a.s. to some $J_* \in L^1(\mathbb{R}_+)$. By combining the convex combination, we can always assume that $(\hat{y}_n, \hat{Y}_n, \hat{J}_n) \in$ conv$\{(y_k, Y_k, \tilde{U}(Y_k)), \, k \geq n\}$.

We now prove that the latter convergence holds in $L^1$ and that $\tilde{U}(Y_*) = J_*$. Because $\tilde{U}$ is convex, we have $\hat{J}_n \geq \tilde{U}(\hat{Y}_n)$ and therefore $J_* \geq \tilde{U}(Y_*)$. On the other hand, it follows from (5.6) and the uniform integrability of $(Y_n)_n$ that $EJ_* = E\tilde{U}(Y_*)$. Hence, $\tilde{U}(Y_*) = J_*$. Finally, since $(\hat{J}_n)_n$ is nonnegative, converges $P$-a.s. to $J_*$ and $E\hat{J}_n \to EJ_*$, the convergence holds in $L^1$.

STEP 4. It follows from (2.8), (5.7) and (5.5) that

$$V(x) \leq W(x) = \lim_{n \to \infty} V_n(x) \leq \lim_{n \to \infty} V_n^s(x) \leq V(x),$$

which concludes the proof. □

We continue the proof of Theorem 3.1 by turning to the sequences $(X_n)_n$ and $(X_n^s)_n$. Set

$$(5.8) \qquad \begin{aligned} Z_n &:= (X_n - B_n)\mathbf{1}_{\{Y_* > 0\}} + L\mathbf{1}_{\{Y_* = 0\}}, \\ Z_n^s &:= (X_n^s - B_n)\mathbf{1}_{\{Y_* > 0\}} + L\mathbf{1}_{\{Y_* = 0\}}, \end{aligned}$$

where $L := \inf\{l \in [0, \infty] : U(l) = U(\infty)\} \in \mathbb{R} \cup \{+\infty\}$. We then use the convention

$$L \times 0 = 0 \qquad \text{so that } \tilde{U}(0) = U(L) - L \times 0 \text{ is valid.}$$

LEMMA 5.5. There is a sequence $(\hat{Z}_n, \hat{Z}_n^s) \in \text{conv}\{(Z_k, Z_k^s), \, k \geq n\}$ such that

$$(\hat{Z}_n, \hat{Z}_n^s) \to (X_* - B, X_*^s - B), \qquad P\text{-a.s.}$$
$$\text{with } E[X_* Y_*] \leq xy_* \quad \text{and} \quad E[X_*^s Y_*] \leq xy_*.$$

Moreover, $X_* - B \leq L$, $X_*^s - B \leq L$ and $X_* = X_*^s = L$ on $\{Y_* = 0\}$.

PROOF.

STEP 1. We first prove the required result for the sequence $(Z_n)_n$. Recall that on the event set $\{Y_* > 0\}$, $Z_n \in -\partial \tilde{U}_n(Y_n) = -\partial \tilde{U}(Y_n)$ for all $n$ [see (3.1)]. By Lemma 2.3 and the convexity of $\tilde{U}$, it follows that for all $Z_* \in -\partial \tilde{U}(Y_*)$,

$$(5.9) \qquad \begin{aligned} Z_n^- Y_* \mathbf{1}_{\{Y_* > 0\}} &\leq |Z_n| Y_n \mathbf{1}_{\{Y_n > Y_* > 0\}} + |Z_*| Y_* \mathbf{1}_{\{Y_* > 0\}} \mathbf{1}_{\{Y_n \leq Y_*\}} \\ &\leq C(\tilde{U}_n(Y_n) + \tilde{U}(Y_*)). \end{aligned}$$

By Lemma 5.4, (5.2) and the fact that $x > \|B\|_\infty$, this provides

$$\sup_n E[Z_n^- Y_*] < \infty.$$

Also notice that the equality $EY_* = y_*$ implies that

$$(5.10) \quad E[Z_n Y_* \mathbf{1}_{\{Y_* > 0\}}] = E[(X_n - B_n) \mathbf{1}_{\{Y_* > 0\}} Y_*] \le x y_* - E[Y_* B]$$

since $X_n \in \mathcal{C}(x_n)$. It follows that $\sup_k EY_* |Z_k| \mathbf{1}_{\{Y_* > 0\}} < \infty$. Hence, there exists a convex combination $Y_* \hat{Z}_n \mathbf{1}_{\{Y_* > 0\}} \in \text{conv}\{Y_* Z_k \mathbf{1}_{\{Y_* > 0\}}, k \ge n\}$ that converges $P$-a.s. It follows that there exist some $Z_* (=: X_* - B)$ such that $\hat{Z}_n \to Z_*$, $P$-a.s., $Z_* \le L$ and $Z_* \mathbf{1}_{\{Y_* = 0\}} = L$. By combining the convex combinations, we may assume that the coefficients that define $\hat{Z}_n$ and $\hat{J}_n$ are the same. Recall from Lemma 5.4 that $\hat{J}_n$ is uniformly integrable. Then, we deduce from (5.9) that

$$(5.11) \qquad \text{the sequence } \hat{Z}_n^- Y_* \mathbf{1}_{\{Y_* > 0\}} \text{ is uniformly integrable.}$$

Since $E[Z_n Y_* \mathbf{1}_{\{Y_* > 0\}}] \le x y_* - E[BY_*]$, it follows from Fatou's lemma that

$$E[Z_* Y_*] = E[Z_* Y_* \mathbf{1}_{\{Y_* > 0\}}] \le x y_* - E[BY_*].$$

STEP 2. Since $(y_*, Y_*) \in \tilde{\mathbf{Y}}_+$ and $X_n^s \in \mathcal{X}_b(x_n)$, we clearly have $EY_* Z_n^s \le y_* x - EY_* B$. We then observe that $(Z_n^s)^- \le Z_n^-$ and the required results of the sequence $(Z_n^s)_n$ follow by the same argument as above. $\quad \square$

LEMMA 5.6. Let $X_* = X_*^s$ and

$$X_* \in B - \partial \tilde{U}(Y_*), \qquad P\text{-a.s.}, \quad EX_* Y_* = x y_*,$$

so that $EU(X_* - B) = V(x) = W(x) = E[\tilde{U}(Y_*) - Y_* B + x y_*]$. Moreover,

$$Y_* \hat{Z}_n^s \to Y_*(X_* - B) \qquad \text{in } L^1(P).$$

PROOF.

STEP 1. We first prove that

$$(5.12) \qquad X_* \in B - \partial \tilde{U}(Y_*), \qquad P\text{-a.s.} \quad \text{and} \quad EX_* Y_* = x y_*.$$

Notice that by (3.1) and Lemma 2.3,

$$
\begin{aligned}
(5.13) \quad U(Z_n)^+ \mathbf{1}_{\{Y_* > 0\}} \mathbf{1}_{\{Y_n > 0\}} &= U_n(Z_n)^+ \mathbf{1}_{\{Y_* > 0\}} \mathbf{1}_{\{Y_n > 0\}} \\
&\le (\tilde{U}_n(Y_n) + |Z_n| Y_n) \mathbf{1}_{\{Y_* > 0\}} \mathbf{1}_{\{Y_n > 0\}} \\
&\le C \tilde{U}_n(Y_n).
\end{aligned}
$$

Let $(\mu_n^k)$ be the coefficients of the convex combination defined in Lemma 5.4(i). Since, by Remark 3.7, $Y_n > 0$ whenever $U(\infty) = \infty$, we deduce from the above inequalities that

$$\left\{ \sum_{k \ge n} \mu_n^k U(Z_k) \right\}^+ \mathbf{1}_{\{Y_* > 0\}} \le C(1 + \hat{J}_n),$$

which is uniformly integrable by Lemma 5.4. It follows from Lemma 5.4, (5.3), the definition of $Z_n$ in (5.8), Fatou's lemma, the concavity of $U$ and Lemma 5.5 that

$$
\begin{aligned}
W(x) &= \lim_{n\to\infty} E\left[\sum_{k\geq n}\mu_n^k U(X_k - B_k)\right] \\
&\leq \lim_{n\to\infty} E\left[\sum_{k\geq n}\mu_n^k (U(Z_k)\mathbf{1}_{\{Y_*>0\}} + U(\infty)\mathbf{1}_{\{Y_*=0\}})\right] \\
(5.14)\qquad &\leq E\left[\limsup_{n\to\infty}\sum_{k\geq n}\mu_n^k (U(Z_k)\mathbf{1}_{\{Y_*>0\}} + U(\infty)\mathbf{1}_{\{Y_*=0\}})\right] \\
&\leq E\left[\lim_{n\to\infty} U(\hat{Z}_n)\right] \\
&= EU(Z_*).
\end{aligned}
$$

By Lemmas 5.5 and 5.4(ii), we get that

$$
W(x) \leq EU(Z_*) \leq E[\tilde{U}(Y_*) + Z_*Y_*] \leq E[\tilde{U}(Y_*) - Y_*B + xy_*] = W(x).
$$

Then equality holds and (5.12) follows.

STEP 2.   From (5.11) and the fact that $Z_n^s \geq Z_n$, we see that

(5.15)      the sequence $\{Y_*(\hat{Z}_n^s)^-,\, n \geq 0\}$ is uniformly integrable.

We also recall that $X_n^s \in \mathcal{X}_b(x_n)$ and therefore

$$
(5.16)\qquad E\hat{Z}_n^s Y_* \leq xy_* - EY_*B.
$$

It then follows from Fatou's lemma together with Step 1 of this proof that $EX_*^s Y_* \leq xy_* = EX_* Y_*$ so that $E(X_*^s - X_*)Y_* \leq 0$. Since $X_*^s - X_* \geq 0$ and $X_* = X_*^s$ on $\{Y_* = 0\}$ by Lemma 5.5, this provides $X_* = X_*^s$, $P$-a.s.

STEP 3.   It remains to prove the $L^1(P)$ convergence of the sequence $(Y_*\hat{Z}_n^s)_n$. To see this, apply Fatou's lemma in (5.16) and use the equality $EX_* Y_* = xy_*$. The result is

$$
E[Y_*\hat{Z}_n^s] \to E[Y_*(X_* - B)].
$$

Since $\hat{Z}_n^s \to Z_*^s = Z_*$, $P$-a.s. by Step 2 of this proof, the required result follows from (5.15).   □

LEMMA 5.7.   *We have*

$$
\sum_{k\geq n}\mu_n^k U(X_n - B_n) \to U(X_* - B) \qquad in\ L^1,
$$

*where $(\mu_n^k)$ are the coefficients of the convex combination defined in Lemma 5.4(i).*

PROOF. Set $I_n = U(X_n - B_n)$. By Remark 3.7, $Y_n > 0$ whenever $U(\infty) = \infty$. From Lemma 2.3 and (5.3), it follows that

$$
\begin{aligned}
[I_n]^+ &\leq [U(X_n - B_n)]^+ \mathbf{1}_{\{Y_n > 0\}} + C \\
&\leq \tilde{U}_n(Y_n) \mathbf{1}_{\{Y_n > 0\}} + |X_n - B_n| Y_n \mathbf{1}_{\{Y_n > 0\}} + C \\
&\leq C(1 + \tilde{U}_n(Y_n))
\end{aligned}
\tag{5.17}
$$

for some constant $C > 0$. Hence, by Lemma 5.2, (5.2) and the fact that $x > \|B\|_\infty$, it follows that $\sup_n E[I_n^+] < \infty$. Since

$$
\sup_n |EI_n| = \sup_n |EU_n(X_n - B_n)| = \sup_n |V_n(x)| < \infty,
$$

by Lemma 5.4, it follows that

$$
\sup_n E|I_n| < \infty.
$$

Hence, we can find a sequence $\hat{I}_n \in \text{conv}\{I_k,\ k \geq n\}$ that converges $P$-a.s. to some $I_*$. By combining the convex combinations, we can assume that the coefficients that define $\hat{I}_n$, $\hat{Z}_n$ and $\hat{J}_n$ are the same. Since by concavity of $U$, $\hat{I}_n \leq U(\hat{Z}_n)$, we have

$$
I_* \leq U(Z_*).
$$

Moreover, because the sequence $(\hat{J}_n)_n$ is uniformly integrable (see Lemma 5.4), it follows from (5.17) that $([\hat{I}_n]^+)_n$ is uniformly integrable. Using (5.14), $W(x) = EU(Z_*)$ (see Lemma 5.6) and Fatou's lemma, we obtain that $EU(Z_*) \leq EI_*$ and therefore

$$
U(Z_*) = I_*.
$$

Since, by (5.14), $E\hat{I}_n \to EU(Z_*)$, we obtain that $\hat{I}_n \to U(Z_*)$ in $L^1$. $\quad\square$

We are now able to complete the proof of Theorem 3.1(ii).

COROLLARY 5.1. Let $X_n^{s'} := X_n^s - n/2$ and $\hat{X}_n^{s'} := \sum_{k \geq n} \mu_n^k X_k^{s'}$, where $(\mu_n^k)$ are the coefficients of the convex combination defined in Lemma 5.4(i). Then

$$
\hat{X}_n^{s'} \in \mathcal{X}_b(x) \quad \text{and} \quad U(\hat{X}_n^{s'} - B) \to U(X_* - B) \qquad \text{in } L^1.
$$

PROOF. By Lemma 5.5, Lemma 5.6 and the concavity of $U$,

$$
\begin{aligned}
\sum_{k \geq n} \mu_n^k U(X_n - B_n) &\leq U(\hat{X}_n^{s'} - B) \\
&= U(\hat{Z}_n^s) \to U(X_*^s - B) = U(X_* - B), \qquad P\text{-a.s.}
\end{aligned}
$$

By Lemma 5.4 and the fact that $\hat{X}_n^{s'} \in \mathcal{X}_b(x)$, this provides

$$V(x) = \lim_n E \sum_{k \geq n} \mu_n^k U(X_n - B_n) \leq \lim_n EU(\hat{X}_n^{s'} - B) \leq V(x).$$

The required result follows from the $L^1(P)$ convergence result of Lemma 5.7. $\square$

Items (ii) and (iii) of Theorem 3.1 are obtained by combining Corollary 5.1 with Lemma 5.6. We conclude the proof of Theorem 3.1 by verifying item (iv).

LEMMA 5.8. *Assume that $Y_* > 0$, $P$-a.s. Then $X_* = X_T^{x,\theta}$ for some $\theta \in L(S)$, where $X^{x,\theta}$ is a uniformly integrable martingale under $Q_* := \frac{Y_*}{y_*} \cdot P$.*

PROOF. Set $H_* = Y_*/y_*$. For $t \leq T$, define

$$M_t := E[H_* X_* | \mathcal{F}_t].$$

Since, $E[H_* | X_*|] < \infty$ by Lemmas 5.6 and 2.3, this defines a process $M$ which is a uniformly integrable martingale under $Q_* := H_* \cdot P$. Also notice from Lemma 5.6 that $M_0 = x$. Finally recall that $\hat{X}_n^{s'} \in \mathcal{X}_b(x)$ and, by Lemma 5.6 and $Y_* > 0$,

$$\hat{X}_n^{s'} = \hat{Z}_n^s + B \to X_* \qquad \text{in } L^1(Q_*).$$

The proof is now completed by the same argument as in Step 10 of [18]. $\square$

The proof of Theorem 3.1 is complete. We conclude this section with the proof of Proposition 3.1.

PROOF OF PROPOSITION 3.1. The case $U(\infty) = \infty$ already was discussed in Section 3. We then assume that $U$ is bounded from above.

STEP 1. We first prove that $\partial \tilde{U}(0) = \{-\infty\}$. To see this, observe that because $U$ is bounded from above, nondecreasing and concave, we have that $\partial U(+\infty) = \{0\}$. Now suppose that $0 \in \partial U(x)$ for some finite $x$. Then $U(x) = U(\infty)$ by concavity of $U$ and $L \leq x$, which contradicts the assumption of the lemma. The required result follows from the classical connection between the generalized gradients of $U$ and $\tilde{U}$.

STEP 2. Let $(y_*, Y_*)$ be the solution of $W(x)$ and define $(y_\varepsilon, Y_\varepsilon) := \varepsilon(\bar{y}, \bar{Y}) + (1 - \varepsilon)(y_*, Y_*)$ for some $\varepsilon \in (0, 1/2)$. By convexity of $\tilde{U}$, we have $(y_\varepsilon, Y_\varepsilon) \in \bar{\mathbf{Y}}_+$ and $\tilde{U}(Y_\varepsilon) \in L^1$. Set

$$X_\varepsilon := \operatorname{ess\,inf}\{X \in L^0 : X \in B - \partial \tilde{U}(Y_\varepsilon)\},$$

and observe that $B - X_\varepsilon \in \partial\tilde{U}(Y_\varepsilon)$ and $X_\varepsilon \to X_0$, $P$-a.s. with $X_0 := \operatorname{ess\,inf}\{X \in L^0 : X \in B - \partial\tilde{U}(Y_*)\}$. We now use the optimality of $(y_*, Y_*)$ together with the convexity of $\tilde{U}$. The result is

$$
\begin{aligned}
(5.18) \quad 0 &\geq \frac{1}{\varepsilon}[E\ (\tilde{U}(Y_*) + y_*x - Y_*B) - E(\tilde{U}(Y_\varepsilon) + y_\varepsilon x - Y_\varepsilon G)] \\
&\geq E(Y_* - \bar{Y})(B - X_\varepsilon) + (y_* - \bar{y})x.
\end{aligned}
$$

We prove later that

$$
(5.19) \quad ([(Y_* - \bar{Y})(B - X_\varepsilon)]^-)_{0 < \varepsilon < 1/2} \text{ is uniformly bounded in } L^1,
$$

$$
E[\bar{Y}[X_0 - B]^-] < \infty
$$

and

$$
(5.20) \quad E[Y_*(X_0 - B)] < \infty,
$$

so that (5.18) implies that $E\bar{Y}(X_0 - B)^+ < \infty$. Since $\bar{Y} > 0$, $P$-a.s. and, by Step 1 of this proof, $X_0 - B = +\infty$ on $\{Y_* = 0\}$ this proves that $Y_* > 0$, $P$-a.s.

STEP 3. We now prove (5.19). Since $Y_\varepsilon > 0$, $\tilde{U}$ is convex and $x \mapsto (x)^-$ is nonincreasing, it follows that for all $Z \in -\partial\tilde{U}(\bar{Y})$ and $Z_* \in -\partial\tilde{U}(Y_*)$,

$$
[\bar{Y}(X_\varepsilon - B)]^- \leq \bar{Y}|Z| + Y_*|Z_*|\mathbf{1}_{Y_* > 0}.
$$

By the same type of argument, we obtain that for all $Z_\varepsilon \in -\partial\tilde{U}((1-\varepsilon)Y_*)$,

$$
[-Y_*(X_\varepsilon - B)]^- \leq Y_*|Z_\varepsilon| \leq 2(1-\varepsilon)Y_*|Z_\varepsilon|.
$$

By Lemma 2.3, this provides

$$
\begin{aligned}
[(\bar{Y} - Y_*)(X_\varepsilon - B)]^- &\leq [\bar{Y}(X_\varepsilon - B)]^- + [-Y_*(X_\varepsilon - B)]^- \\
&\leq C\tilde{U}(\bar{Y}) + C\tilde{U}(Y_*) + 2C\tilde{U}((1-\varepsilon)Y_*)\mathbf{1}_{Y_* > 0} \\
&\leq C\tilde{U}(\bar{Y}) + C\tilde{U}(Y_*) + 2C^2\tilde{U}(Y_*) \in L^1.
\end{aligned}
$$

The previous inequalities also prove the second claim of (5.19) since $X_\varepsilon \to X_0$, $P$-a.s.

STEP 4. It remains to prove (5.20). Since $X_0$ is valued in $B - \partial\tilde{U}(Y_*)$ and $X_0 \leq X_*$, it follows from the definition of $\tilde{U}$ together with the nondecrease of $U$ that

$$
\tilde{U}(Y_*) = U(X_0 - B) - Y_*(X_0 - B) \leq U(X_* - B) - Y_*(X_0 - B),
$$

so that $EY_*(X_0 - B) \leq V(x) - E\tilde{U}(Y_*) < \infty$. $\quad\square$

**6. Utility functions with bounded negative domain.** In this section, we proceed to the proof of Theorem 3.2 which was the starting point of the proof of Theorem 3.1. We warn the reader that many notations from the previous sections will be used in this section for different objects.

The effective domains of the utility function and the associated Fenchel transform are now assumed to satisfy

$$\mathrm{cl}(\mathrm{dom}(U)) = [-2\beta, \infty) \quad \text{and} \quad \mathrm{cl}(\mathrm{dom}(\tilde{U})) = \mathbb{R}_+.$$

Recall that we have assumed

$$(6.1) \qquad\qquad\qquad U(+\infty) > 0,$$

so that $\tilde{U}(0+) > 0$. The following remark collects some properties of $\tilde{U}$.

REMARK 6.1. (i) The function $y \mapsto \tilde{U}(y) - 2\beta y$ is nonincreasing and positive near 0.

(ii) By Lemma 2.1,

$$\mathrm{AE}_0(\tilde{U}(\cdot)) < \infty \quad \implies \quad \mathrm{AE}_0(\tilde{U}(\cdot) - 2\beta\cdot) < \infty.$$

It follows from Lemma 4.1 in [5] that the asymptotic elasticity condition $\mathrm{AE}_0(\tilde{U}) < \infty$ is equivalent to the existence of two constants $\gamma > 0$ and $y_0 > 0$ such that

$$\tilde{U}(\mu y) - 2\beta\mu y \leq \mu^{-\gamma}(\tilde{U}(y) - 2\beta y) \qquad \text{for all } \mu \in (0,1] \text{ and } y \in (0, y_0].$$

(iii) Applying the latter characterization to $y_0$ and using the nonincrease property (i), we see that

$$\tilde{U}(y_0) - 2\beta y_0 \leq \tilde{U}(\mu y_0) - 2\beta\mu y_0 \leq \mu^{-\gamma}[\tilde{U}(y_0) - 2\beta y_0]$$

for any arbitrary $\mu \in (0,1)$. This proves that $\tilde{U}(y_0) - 2\beta y_0 \geq 0$ and, by (i),

$$\tilde{U}(y) - 2\beta y \geq 0 \qquad \text{for all } y \in (0, y_0].$$

(iv) Fix $\bar{y} \in (0, \infty)$. Then, using a compactness argument, we deduce from the characterization (ii) of the asymptotic elasticity condition $\mathrm{AE}_0(\tilde{U}) < \infty$ that there exist positive constants $\gamma > 0$ and $C_{\bar{y}} > 0$ such that

$$\tilde{U}(\mu y) - 2\beta\mu y \leq \mu^{-\gamma}[\tilde{U}(y) - 2\beta y + C_{\bar{y}}] \qquad \text{for all } \mu \in [1/2, 1] \text{ and } y \in (0, \bar{y}).$$

6.1. *Approximation by quadratic inf convolution.* The main difficulty arises from the nonsmoothness of $\tilde{U}$ inherited from $U$. To handle this problem, we introduce the quadratic inf convolution:

$$\tilde{U}_n(y) := \beta y + \inf_{z \geq 0}\left(\tilde{U}(z) - \beta z + \frac{n}{2}|y - z|^2\right).$$

Then $\tilde{U}_n$ is finitely defined on $\mathbb{R}$, strictly convex and

$$(6.2) \qquad \tilde{U}_n(y) \le \tilde{U}(y) \qquad \text{for all } y \ge 0.$$

We report from [5] the following properties of $\tilde{U}_n$ which will be used in the subsequent analysis.

PROPERTY 6.1. *For all $y \in \mathbb{R}$, there exists a unique $z_n(y) \ge 0$ such that*

$$\tilde{U}_n(y) = \tilde{U}(z_n(y)) - \beta(z_n(y) - y) + \frac{n}{2}|z_n(y) - y|^2.$$

PROPERTY 6.2. (i) *For all $x > 0$ and $y \in \mathbb{R}$, we have*

$$|z_n(y) - y|^2 \le \frac{4}{n}[\tilde{U}_n(y) - \beta y + xy + C]$$

*for some constant $C$.*

(ii) *Let $(y_n)_n$ be a sequence converging to $y \in \mathrm{dom}(\tilde{U})$. Then*

$$z_n(y_n) \to y.$$

(iii) *Let $(y_n)_n$ be a sequence converging to $y$. Suppose further that $z_n(y_n) \to y$. Then*

$$\tilde{U}_n(y_n) \to \tilde{U}(y).$$

PROPERTY 6.3. *Function $\tilde{U}_n$ is continuously differentiable on $\mathbb{R}$ and*

$$D\tilde{U}_n(y) = n(y - z_n(y)) + \beta \in \partial \tilde{U}(z_n(y)).$$

REMARK 6.2. From Remark 6.1 and Property 6.3 of the inf convolution, we deduce that

$$y \mapsto \tilde{U}_n(y) - 2\beta y \qquad \text{is nonincreasing.}$$

PROPERTY 6.4. *Suppose that $\mathrm{AE}_0(\tilde{U}) < \infty$. Then there exist some $y_0 > 0$ and some positive constants $\gamma$ and $C$ such that, for all $n \ge 1$,*

$$\tilde{U}_n(\mu y) - \beta \mu y \le \mu^{-\gamma}[C + \tilde{U}_n(y) - \beta y] \qquad \text{for all } \mu \in [1/2, 1] \text{ and } y \in (0, y_0]$$

*and*

$$-(D\tilde{U}_n(y) - \beta)y \le C(1 + \tilde{U}_n(y) - \beta y) \qquad \text{for all } y \in (0, y_0].$$

PROOF. The second inequality follows from the first by the same type of arguments as in the proof of Lemma 2.3(ii) (see the Appendix). We now concentrate on the first inequality. We set $g_n(y) := \tilde{U}_n(y) - \beta y$ and $g(y) := \tilde{U}(y) - \beta y$.

STEP 1. Let $y_0 > 0$ be defined as in Remark 6.1. Fix $0 < y \leq y_0$ and define

$$f_n(z) := g(z) + \frac{n}{2}|y - z|^2.$$

We first prove that $f_n$ is increasing on $(z_0, \infty)$, where

$$z_0 := 2y_0 + |\beta - \max \partial g(y_0)| < \infty$$

is independent of $n \geq 1$ and $0 < y \leq y_0$. Consider some arbitrary $z \geq z_0$ and $q_1 \in \partial f_n(z)$. Then there exist some $q_2 \in \partial(g - \beta \cdot)(z)$ such that

$$\frac{1}{n}q_1 = \frac{1}{n}(q_2 + \beta) + (z - y).$$

Since the map $(g - \beta \cdot)$ is nonincreasing, by Remark 6.1(i), it follows that $q_2 \leq 0$. Since it is also convex and $z \geq y_0 \geq y$, we get

$$\frac{1}{n}q_1 \geq q_2 + (z - y) \geq \max \partial g(y_0) - \beta + z - y_0 \geq z - z_0 + y_0 > 0$$

since $y_0 > 0$ and $z \geq z_0$. This proves that, for all $n \geq 1$ and $0 < y \leq y_0$, $f_n$ is increasing on $[z_0, \infty)$ and therefore

$$\tag{6.3} g_n(\mu y) = \inf_{0 \leq z \leq z_0} \left( g(z) + \frac{n}{2}|\mu y - z|^2 \right)$$

$$\text{for all } (y, \mu) \in (0, y_0] \times [1/2, 1], \ n \geq 1.$$

STEP 2. Fix $(y, \mu) \in (0, y_0] \times [1/2, 1]$. By (6.3), we see that

$$g_n(\mu y) = \inf_{0 \leq z \leq z_0} \left( g(z) + \frac{n}{2}|\mu y - z|^2 \right) = \inf_{0 \leq z \leq 2z_0} \left( g(\mu z) + \mu^2 \frac{n}{2}|y - z|^2 \right),$$

where the second equality is obtained by a trivial change of variable and the fact that $\mu \geq 1/2$. Using Remark 6.1(iv) with $\bar{y} = 2z_0$, we deduce that there exist some $C > 0$ and $\gamma > 0$, such that

$$g_n(\mu y) \leq \inf_{0 \leq z \leq 2z_0} \left( \mu^{-\gamma}(C + g(z) - \beta z) + \beta \mu z + \mu^2 \frac{n}{2}|y - z|^2 \right).$$

Since $\mu \leq \mu^{-\gamma}$ and $\mu^{2+\gamma} \leq 1$, this provides

$$g_n(\mu y) \leq \inf_{0 \leq z \leq 2z_0} \left( \mu^{-\gamma}(C + g(z)) + \mu^2 \frac{n}{2}|y - z|^2 \right)$$

$$\leq \mu^{-\gamma} \left[ C + \inf_{0 \leq z \leq 2z_0} \left( g(z) + \frac{n}{2}|y - z|^2 \right) \right] = \mu^{-\gamma}(C + g_n(y)),$$

where the last inequality follows from (6.3) again. □

By substituting $\tilde{U}_n$ for $\tilde{U}$ in the definition of the dual problem

$$(6.4) \qquad W(x) := \inf_{(y,Y) \in \mathbf{Y}_+} E[\tilde{U}(Y) - YB + xy]$$

of Theorem 3.2, we define a sequence of approximate dual problems:

$$(6.5) \qquad W_n(x) := \inf_{(y,Y) \in \mathbf{Y}_+} E[\tilde{U}_n(Y) - YB + xy].$$

6.2. *Existence in the dual problem.* The purpose of this section is to prove that the approximate dual problem $W_n(x)$ has a solution, for each $n$, and to define a solution for the dual problem $W(x)$ as a limit of these solutions in some appropriate sense.

The following preliminary result will be used frequently.

LEMMA 6.1. *Let $\beta = 0$. Then there exists a sequence of functions $(\phi_n)_{1 \le n \le \infty}$ such that, for all sufficiently large $n$,*

$$\phi_n : (-\tilde{U}_n(0), +\infty) \to (0, \infty) \qquad with \lim_{x \to \infty} \frac{\phi_n(x)}{x} = \infty$$

*and*

$$E[\phi_n(\tilde{U}_n(Y)^-)] \le C + y \qquad for\ all\ (y, Y) \in \mathbf{Y}_+\ with\ y > 0$$

*for some $C > 0$ independent of $n$, with the convention $\tilde{U}_\infty = \tilde{U}$. In particular, for all $M > 0$ and large $n$, the family $\{\tilde{U}_n(Y)^-, (y, Y) \in \mathbf{Y}_+, |y| \le M\}$ is uniformly integrable.*

The proof of this result is reported in Section 6.4. We now establish existence in the approximate dual problems $W_n$ and convergence of these solutions (in some sense) to some solution of $W(x)$. These results will be established under the following assumptions.

STANDING ASSUMPTIONS OF SECTION 6.2. $\|B\|_\infty \le \beta, x > 0$ and $W(x) < \infty$.

LEMMA 6.2. *For sufficiently large $n$, existence holds for the problem $W_n(x)$, that is,*

$$W_n(x) = E[\tilde{U}_n(Y_n) + y_n x - Y_n B] \qquad for\ some\ (y_n, Y_n) \in \mathbf{Y}_+.$$

PROOF. Let $n \ge 1$ be a fixed integer and let $(y_k, Y_k)_k$ be a minimizing sequence of $W_n(x)$. Then, from (6.2), we have

$$(6.6) \quad \begin{aligned} &-E[\tilde{U}_n(Y_k) - Y_k B]^- + xy_k \\ &\qquad \le E[\tilde{U}_n(Y_k) - Y_k B] + xy_k \le W_n(x) + 1 \le W(x) + 1. \end{aligned}$$

STEP 1. We first prove that the sequence $(y_k)_k$ is bounded so that, by Remark 5.1, there is a sequence $(\hat{y}_k, \hat{Y}_k) \in \text{conv}\{(y_j, Y_j), j \geq k\}$ which converges $P$-a.s. to some $(\hat{y}, \hat{Y}) \in \mathbf{Y}_+$.

(i) The case $\beta > 0$ is easily dealt with since, with the notation of Property 6.1,

$$(6.7) \qquad \tilde{U}_n(Y_k) \geq \tilde{U}(z_n(Y_k)) - \beta z_n(Y_k) + \beta Y_k \geq U(-\beta) + \beta Y_k,$$

so that (6.6) together with the condition $\|B\|_\infty \leq \beta$ provide

$$xy_k \leq U(-\beta)^- + W(x) + 1.$$

Since $x$ is positive and $y_k$ is nonnegative, this proves that the sequence $(y_k)_k$ is bounded.

(ii) We then concentrate on the case $\beta = 0$. Let $\phi_n$ be the function introduced in Lemma 6.1. Then for all $\varepsilon > 0$, there exists some $x_0 > 0$ such that

$$\frac{\phi_n(x)}{x} \geq \frac{1}{\varepsilon} \qquad \text{for } x \geq x_0,$$

and then,

$$x \leq x_0 + \varepsilon \phi_n(x)\mathbf{1}_{\{x \geq x_0\}} \leq x_0 + \varepsilon \phi_n(x) \qquad \forall x \geq 0,$$

for sufficiently large $x_0$ and $n$. Using Lemma 6.1, we then compute that, for some $C > 0$,

$$E\tilde{U}_n(Y_k)^- \leq x_0 + \varepsilon E\phi_n(\tilde{U}_n(Y_k)^-) \leq x_0 + \varepsilon(C + y_k).$$

Plugging this inequality in (6.6), we obtain

$$(6.8) \qquad (x - \varepsilon)y_k \leq W(x) + 1 + x_0 + \varepsilon C.$$

By choosing $\varepsilon = x/2 > 0$, we see that the sequence $(y_k)_k$ is bounded.

STEP 2. Combining Lemma 6.1, (6.7) and $\beta \geq \|B\|_\infty$, we see that the sequence $\{(\tilde{U}_n(\hat{Y}_k) - \hat{Y}_k B)^-, k \geq 0\}$ is uniformly integrable. Let $(\mu_k^j)$ be the coefficients of the convex combination defining $(\hat{Y}_k)$. By Fatou's lemma, together with the convexity of $\tilde{U}_n$, we get

$$\begin{aligned} W_n(x) &\leq E[\tilde{U}_n(\hat{Y}) - \hat{Y}B] + x\hat{y} \\ &\leq \liminf_{k \to \infty} E[\tilde{U}_n(\hat{Y}_k) - \hat{Y}_k B] + x\hat{y}_k \\ &\leq \liminf_{k \to \infty} \sum_{j \geq k} \mu_k^j E\tilde{U}_n(Y_j) - Y_j B + xy_j = W_n(x), \end{aligned}$$

since $(y_j, Y_j)_j$ is a minimizing sequence of $W_n(x)$. This proves that $(\hat{y}, \hat{Y})$ is a solution of $W_n(x)$. $\square$

REMARK 6.3. For later use, observe that the same arguments as in Step 2 of the above proof show that, for sufficiently large $n$,

the family $\{(\tilde{U}_n(Y) - YB)^- : (y, Y) \in \mathbf{Y}_+, |y| \leq M\}$ is uniformly integrable for all $M > 0$.

The next lemma completes the proof of Theorem 3.2(i).

LEMMA 6.3. Let $(y_n, Y_n)$ be a solution of $W_n(x)$. Then there exists a sequence $(\bar{y}_n, \bar{Y}_n) \in \text{conv}((y_k, Y_k), k \geq n)$ such that

$$(6.9) \qquad (\bar{y}_n, \bar{Y}_n) \to (y_*, Y_*) \in \mathbf{Y}_+, \qquad P\text{-}a.s.$$

Moreover, $(y_*, Y_*)$ is a solution of the problem $W(x)$.

PROOF.

STEP 1. We first argue as in the previous proof to show that the sequence $(y_n)_n$ is bounded so that, by Remark 5.1, there is a sequence $(\bar{y}_n, \bar{Y}_n) \in \text{conv}\{(y_j, Y_j), j \geq n\}$ which converges $P$-a.s. to some $(y_*, Y_*) \in \mathbf{Y}_+$.

By definition of $(y_n, Y_n)$, we have

$$-E[\tilde{U}_n(Y_n) - Y_n B]^- + xy_n \leq E[\tilde{U}_n(Y_n) - Y_n B] + xy_n = W_n(x) \leq W(x).$$

The case $\beta > 0$ is easily solved by observing that $\tilde{U}_n(Y_n) - Y_n B \geq U(-\beta)$ as in (6.7). As for the case $\beta = 0$, we again argue as in the previous proof to derive the analogue of (6.8) with $\varepsilon = x/2 > 0$:

$$xy_n \leq 2W(x) + 2x_0 + xC \leq 2W(x) + 2x_0 + xC \qquad \text{for all large } n$$

and some $C > 0$ independent of $n$. This provides the required bound on $(y_n)_n$.

STEP 2. Set $g(y) := \tilde{U}(y) - \beta y$. Using Property 6.1 of the quadratic inf convolution, we see that

$$g(z_n(\bar{Y}_n)) - \bar{Y}_n(B - \beta) = \tilde{U}_n(\bar{Y}_n) - \bar{Y}_n B - \frac{n}{2}|z_n(\bar{Y}_n) - \bar{Y}_n|^2 \leq \tilde{U}_n(\bar{Y}_n) - \bar{Y}_n B.$$

Let $(\lambda_n^j)_{j \geq n}$ be coefficients of the above convex combination that define $(\bar{y}_n, \bar{Y}_n)$ from $(y_j, Y_j)_{j \geq n}$. From the convexity of $\tilde{U}_n$ and the increase of $\tilde{U}_n$ in $n$, we get from the previous inequality

$$(6.10) \quad g(z_n(\bar{Y}_n)) - \bar{Y}_n(B - \beta) \leq \tilde{U}_n(\bar{Y}_n) - \bar{Y}_n B \leq \sum_{j \geq n} \lambda_n^j [\tilde{U}_j(Y_j) - Y_j B].$$

Then, taking expected values, we see that

$$\begin{aligned}
E[g(z_n(\bar{Y}_n)) - \bar{Y}_n(B - \beta)] &\le E[\tilde{U}_n(\bar{Y}_n) - \bar{Y}_n B] \\
&\le \sum_{j \ge n} \lambda_n^j [W_j(x) - x y_j] \\
&\le W(x) - x \bar{y}_n.
\end{aligned}$$

(6.11)

We now use the claim (the proof of which will be carried out in Step 3 below)

(6.12) the sequence $([g(z_n(\bar{Y}_n)) - \bar{Y}_n(B - \beta)]^-)_n$ is uniformly integrable.

Recalling that $g(\cdot) + \beta \cdot = \tilde{U}(\cdot)$ and using Property 6.2 of the quadratic inf convolution, it follows from Fatou's lemma and (6.11) that

$$\begin{aligned}
E[\tilde{U}(Y_*) - Y_* B] + x y_* &\le \liminf_{n \to \infty} E[\tilde{U}_n(\bar{Y}_n) - \bar{Y}_n B] + x \bar{y}_n \\
&\le \limsup_{n \to \infty} E[\tilde{U}_n(\bar{Y}_n) - \bar{Y}_n B] + x \bar{y}_n \le W(x).
\end{aligned}$$

(6.13)

Since $(y_*, Y_*) \in \mathbf{Y}_+$, this proves that $(y_*, Y_*)$ is the solution of the problem $W(x)$.

STEP 3. To complete the proof, it remains to check (6.12). As in the previous proof, the case $\beta > 0$ is easily solved by observing that $g(z_n(\bar{Y}_n)) = \tilde{U}(z_n(\bar{Y}_n)) - \beta z_n(\bar{Y}_n) \ge U(-\beta)$, so that

$$g(z_n(\bar{Y}_n)) - \bar{Y}_n(B - \beta) \ge U(-\beta) + \bar{Y}_n(\beta - B) \ge U(-\beta),$$

since $\|B\|_\infty \le \beta$. We then concentrate on the case $\beta = B = 0$. Let $\phi := \phi_\infty$ be the function introduced in Lemma 6.1. Then

(6.14) $\quad E[\phi(\tilde{U}(z_n(\bar{Y}_n))^-)] \le C + E[z_n(\bar{Y}_n)] \le C + \bar{y}_n + E[z_n(\bar{Y}_n) - \bar{Y}_n].$

By the first part of this proof, the sequence $(\bar{y}_n)_n$ is bounded. We next use Property 6.2(i) of the quadratic inf convolution together with (6.11) and $\beta \ge \|B\|_\infty$ to see that

$$\begin{aligned}
E|z_n(\bar{Y}_n) - \bar{Y}_n|^2 &\le \frac{4}{n} E[C + \tilde{U}_n(\bar{Y}_n) - \beta \bar{Y}_n + x \bar{Y}_n] \\
&\le \frac{4}{n} E[C + \tilde{U}_n(\bar{Y}_n) - B \bar{Y}_n + x \bar{y}_n] \\
&\le \frac{4}{n} [C + W(x)].
\end{aligned}$$

In particular, this proves that the sequence $(E[z_n(\bar{Y}_n) - \bar{Y}_n])_n$ is bounded. Hence the right-hand side term of (6.14) is bounded. Since $\phi(x)/x \to \infty$ as $x \to \infty$, this proves (6.12) by the la Vallée–Poussin theorem. $\quad\square$

REMARK 6.4. For later use, observe that the arguments of Step 4 of the above proof also hold if we replace $(\bar{y}_n, \bar{Y}_n)$ with $(y_n, Y_n)$. It follows that the sequence $([g(z_n(Y_n)) - Y_n(B - \beta)]^-)_n$ is uniformly integrable. Using Property 6.1 as in Step 2, we see that

$$g(z_n(Y_n)) - Y_n(B - \beta) \leq \tilde{U}_n(Y_n) - Y_n B,$$

so that

(6.15)    the sequence $([\tilde{U}_n(Y_n) - Y_n B]^-)_n$ is uniformly integrable.

COROLLARY 6.1. $W_n(x) \to W(x)$.

PROOF. Recall that the sequence $(W_n(x))_n$ is nondecreasing. Since $W_n(x) \leq W(x)$, we have $W_n(x) \to W_\infty(x)$ for some $W_\infty(x) \leq W(x)$. The result is then obtained by combining (6.11) and (6.13) in the above proof.   □

COROLLARY 6.2. Let $(y_n, Y_n)$ be a solution of $W_n(x)$ and let $(y_*, Y_*)$ be the limit defined in Lemma 6.3. Set $J_n := \tilde{U}_n(Y_n) - Y_n B$. Then there exists a sequence $(\hat{y}_n, \hat{Y}_n, \hat{J}_n) \in \text{conv}((y_k, Y_k, J_k), k \geq n)$ such that

$$(\hat{y}_n, \hat{Y}_n) \to (y_*, Y_*), \qquad P\text{-a.s.} \quad and \quad \hat{J}_n \to \tilde{U}(Y_*) - Y_* B \qquad in \ L^1(P).$$

PROOF. From Lemma 6.3, there exists a sequence $(\bar{y}_n, \bar{Y}_n) \in \text{conv}((y_k, Y_k), k \geq n)$ which converges $P$-a.s. to a solution $(y_*, Y_*)$ of $W(x)$. Denote by $(\lambda_k^n, k \geq n)$ the coefficients that define the convex combination and let $\bar{J}_n := \sum_{k \geq n} \lambda_k^n J_k$.

STEP 1. We first prove the existence of a sequence $(\hat{y}_n, \hat{Y}_n, \hat{J}_n) \in \text{conv}((y_k, Y_k, J_k), k \geq n)$ and a random variable $J_* \in L^1(P)$ such that

(6.16)
$$(\hat{y}_n, \hat{Y}_n, \hat{J}_n) \to (y_*, Y_*, J_*) \quad \text{and} \quad E\hat{J}_n \to E\tilde{U}(Y_*) - Y_* B, \qquad P\text{-a.s.}$$

To see this, observe that

$$E\bar{J}_n = \sum_{k \geq n} \lambda_k^n [W_k(x) - xy_*] \quad \Longrightarrow \quad W(x) - xy_* = E\tilde{U}(Y_*) - Y_* B$$

by Corollary 6.1. Also, it follows from (6.10) that

$$\bar{J}_n^- \leq [g(z_n(\bar{Y}_n)) - \bar{Y}_n(B - \beta)]^-,$$

where $g(\cdot) = \tilde{U}(\cdot) - \beta \cdot$. Since the sequence on the right-hand side is uniformly integrable by (6.12), this shows that

(6.17)                    $(\bar{J}_n^-)_n$ is uniformly integrable

and therefore bounded in $L^1$.

Since $|\bar{J}_n| = \bar{J}_n + 2\bar{J}_n^-$, the above arguments show that the sequence $(\bar{J}_n)_n$ is bounded in $L^1$, and (6.16) follows from the Komlòs lemma.

STEP 2. We now prove that

$$(6.18) \qquad J_* = \tilde{U}(Y_*) - Y_* B.$$

By convexity of $\tilde{U}_n$ and increase of $(\tilde{U}_n)_n$, we see that $\hat{J}_n \geq \tilde{U}_n(\hat{Y}_n) - \hat{Y}_n B$. This proves, first, that $[\hat{J}_n]^- \leq [\tilde{U}_n(\hat{Y}_n) - \hat{Y}_n B]^-$ is uniformly integrable by Remark 6.3 and, therefore,

$$EJ_* \leq \lim_{n \to \infty} E\hat{J}_n = E\tilde{U}(Y_*) - Y_* B$$

by Fatou's lemma. This also proves that $J_* \geq \tilde{U}(Y_*) - Y_* B$ by Property 6.2, and (6.18) follows.

STEP 3. In the previous steps, we have proved that $\hat{J}_n \to \tilde{U}(Y_*) - Y_* B$, $P$-a.s, $E\hat{J}_n \to E\tilde{U}(Y_*) - Y_* B$, and $([\hat{J}_n]^-)_n$ is uniformly integrable. This provides that $\hat{J}_n \to \tilde{U}(Y_*) - Y_* B$ in $L^1$. $\square$

6.3. *Existence for the initial problem.* We now turn to the solution of the initial problem $V(x)$. To do this this, we appeal to the following assumptions:

STANDING ASSUMPTIONS OF SECTION 6.3. $\mathcal{M}^e(S) \neq \varnothing$ and $\mathrm{AE}_0(\tilde{U}) < \infty$.

We first start by a characterization of the optimality of $(y_n, Y_n)$ for the problem $W_n(x)$. Recall that $\tilde{U}_n$ is continuously differentiable by Property 6.3.

LEMMA 6.4. *Let $(y_n, Y_n)$ be a solution of $W_n(x)$ and set $X_n := -D\tilde{U}_n(Y_n) + B$. Then:*

(i) $EX_n Y - xy \leq EX_n Y_n - xy_n = 0$ *for all $(y, Y) \in \mathbf{Y}_+$.*

(ii) *There exists a sequence $\hat{X}_n \in \mathrm{conv}(X_k, k \geq n)$ such that $\hat{X}_n \to X_*$ for some $X_*$ in $\mathcal{C}(x)$.*

PROOF.

STEP 1. We first show that (ii) follows easily from (i). Let $Q := Y \cdot P$ be an arbitrary measure in $\mathcal{M}^e(S)$ so that $(1, Y) \in \mathbf{Y}_+$. Since $-3\beta \leq X_n$, we have $E[|X_n|Y] \leq E[X_n Y] + 2E[Y X_n^-] \leq x + 6\beta$ by (5.2) and (i). It follows that the sequence $(X_n)_n$ is bounded in $L^1(Q)$, and the existence of a converging convex combination follows from the Komlòs lemma. Using again (i), we have $E\hat{X}_n Y - xy \leq 0$ for all $(y, Y) \in \mathbf{Y}_+$ and, therefore, $EX_* Y \leq xy$ follows from Fatou's lemma. Clearly, $X_* \geq -3\beta$ and, therefore, $X_* \in \mathcal{C}(x)$.

STEP 2. We prove in Step 3 of this proof that

$$EX_n(Y - Y_n) \leq x(y - y_n) \qquad \text{for all } (y, Y) \in \mathbf{Y}_+.$$

Applying this inequality to $(y, Y) = 2(y_n, Y_n) \in \mathbf{Y}_+$, we see that $EX_nY_n \leq xy_n$. Similarly, by taking $(y, Y) = 2^{-1}(y_n, Y_n) \in \mathbf{Y}_+$, we obtain the converse inequality and then $EX_nY_n = xy_n$. This provides the required result.

STEP 3. Let $(y, Y) \in \mathbf{Y}_+$ be fixed and define for small $\varepsilon > 0$,

$$(y_n^\varepsilon, Y_n^\varepsilon) = (1 - \varepsilon)(y_n, Y_n) + \varepsilon(y, Y) \quad \text{and} \quad X_n^\varepsilon := -D\tilde{U}_n(Y_n^\varepsilon) + B.$$

Clearly, $(y_n^\varepsilon, Y_n^\varepsilon) \in \mathbf{Y}_+$ and as $\varepsilon \searrow 0$, we have $(Y_n^\varepsilon, X_n^\varepsilon) \to (Y_n, X_n)$, $P$-a.s. By the optimality of $(y_n, Y_n)$ for the problem $W_n(x)$ and the convexity of $\tilde{U}_n$, we have

$$0 \geq \varepsilon^{-1}E[\tilde{U}_n(Y_n) - \tilde{U}_n(Y_n^\varepsilon) - B(Y_n - Y_n^\varepsilon)] + \varepsilon^{-1}x(y_n - y_n^\varepsilon)$$
$$\geq EX_n^\varepsilon(Y - Y_n) - x(y - y_n).$$

In the rest of this proof, we show that

$$(6.19) \qquad \text{the sequence } ([X_n^\varepsilon(Y - Y_n)]^-)_\varepsilon \text{ is uniformly integrable,}$$

which provides the required result by sending $\varepsilon$ to zero in the last inequality and using Fatou's lemma.

Let $\alpha$ be a given parameter in $(0, 1/4)$ and $0 \leq \varepsilon \leq \alpha$. We denote $\alpha_\varepsilon := \alpha + \varepsilon$. By convexity of $\tilde{U}_n$ together with Remark 6.2, we see that

$$\tilde{U}_n((1 - \alpha_\varepsilon)Y_n) \geq \tilde{U}_n(Y_n^\varepsilon + \alpha(Y - Y_n)) - \alpha_\varepsilon Y D\tilde{U}_n(Y_n^\varepsilon + \alpha(Y - Y_n))$$
$$\geq \tilde{U}_n(Y_n^\varepsilon + \alpha(Y - Y_n)) - 2\alpha_\varepsilon \beta Y.$$

Using again the convexity of $\tilde{U}_n$, we get

$$(6.20) \qquad \tilde{U}_n((1 - \alpha_\varepsilon)Y_n) \geq \tilde{U}_n(Y_n^\varepsilon) + \alpha D\tilde{U}_n(Y_n^\varepsilon)(Y - Y_n) - 2\beta\alpha_\varepsilon Y$$
$$= J_n^\varepsilon - \alpha X_n^\varepsilon(Y - Y_n) - \alpha_\varepsilon Y(2\beta - B) + (1 - \alpha_\varepsilon)Y_n B,$$

where we set $J_n^\varepsilon := \tilde{U}_n(Y_n^\varepsilon) - Y_n^\varepsilon B$. We now use the asymptotic elasticity condition $AE_0(\tilde{U}) < \infty$ together with Property 6.4 and Remark 6.2 to obtain

$$\tilde{U}_n((1 - \alpha_\varepsilon)Y_n) \leq (1 - \alpha_\varepsilon)^{-\gamma}[C + \tilde{U}_n(Y_n) - \beta Y_n]\mathbf{1}_{\{Y_n \leq y_0\}}$$
$$+ \{\tilde{U}_n((1 - \alpha_\varepsilon)Y_n) - 2(1 - \alpha_\varepsilon)Y_n\beta\}\mathbf{1}_{\{Y_n \geq y_0\}}$$
$$+ (1 - \alpha_\varepsilon)Y_n\beta(1 + \mathbf{1}_{\{Y_n \geq y_0\}})$$
$$\leq C + (1 - \alpha_\varepsilon)^{-\gamma}\tilde{U}_n(Y_n)^+ + 2(1 - \alpha_\varepsilon)Y_n\beta$$

for some $C > 0$. It follows from (6.20) that

$$
\begin{aligned}
\alpha X_n^\varepsilon (Y - Y_n) &\geq J_n^\varepsilon - \alpha_\varepsilon Y(2\beta - B) + (1 - \alpha_\varepsilon) Y_n B \\
&\quad - C - (1 - \alpha_\varepsilon)^{-\gamma} \tilde{U}_n(Y_n)^+ - 2(1 - \alpha_\varepsilon) Y_n \beta \\
&\geq -C - [J_n^\varepsilon]^- - (1 - 2\alpha)^{-\gamma} \tilde{U}_n(Y_n)^+ - 2\alpha Y(2\beta - B) \\
&\quad + (1 - 2\alpha) Y_n(\beta + B) - 3 Y_n \beta,
\end{aligned}
$$

where we used the assumption $\|B\|_\infty \leq \beta$. This provides (6.19) by observing that $Y$, $Y_n$ and $\tilde{U}_n(Y_n)^+$ are integrable, $B$ is bounded, and the family $([J_n^\varepsilon]^-)_\varepsilon$ is uniformly integrable by Remark 6.3. $\quad \square$

LEMMA 6.5. *Let $X_*$ be as in the previous lemma. Then*

$$
E X_* Y_* = x y_*, \qquad X_* \in B - \partial \tilde{U}(Y_*), \qquad P\text{-}a.s. \quad and \quad E U(X_* - B) = V(x).
$$

*Moreover, $V(x) = W(x)$.*

PROOF. Let $(\hat{y}_n, \hat{Y}_n, \hat{X}_n, \hat{J}_n) \in \mathrm{conv}\{(y_k, Y_k, X_k, J_k), k \geq n\}$ be the sequence defined in Corollary 6.2 and Lemma 6.4 (clearly, we can assume that the convex combinations are the same in both results). Define $U_n(x) := \inf_{y \geq 0} \tilde{U}_n(y) + xy$ and observe that $U_n \leq U$. Set

$$
I_n := U_n(X_n - B)
$$

and let $\hat{I}_n$ be the corresponding convex combination.

STEP 1. We claim that

(6.21)             the sequence $(\hat{I}_n^+)_n$ is uniformly integrable.

Before proving this, let us complete the proof of the lemma by repeating the argument of the proof of Lemma 5.6. By Lemma 6.4, $E X_n Y_n = x y_n$ and, therefore,

$$
W_n(x) = E[\tilde{U}_n(Y_n) + x y_n - Y_n B] = E U_n(X_n - B) = E I_n.
$$

Since $W_n(x) \to W(x) = E[\tilde{U}(Y_*) + x y_* - Y_* B]$, it follows from (6.21), Fatou's lemma and the fact that $X_* \in \mathcal{C}(x)$ [see Lemma 6.4(ii)] that

$$
\begin{aligned}
W(x) &= E[\tilde{U}(Y_*) + x y_* - Y_* B] \\
&\leq E\left[ \limsup_n \hat{I}_n \right] \\
&\leq E\left[ \limsup_n U(\hat{X}_n - B) \right] \\
&= E U(X_* - B) \leq V(x) \leq W(x),
\end{aligned}
$$

where we used the fact that $U_n \leq U$ and the concavity of $U$. Then equality holds in the above inequalities and the required results follow.

STEP 2. We now prove (6.21). We first need a preliminary result. Fix $\varepsilon > 0$ and observe that

$$U_n(x) \leq \tilde{U}_n(\varepsilon) + \varepsilon x \qquad \text{for all } x > -2\beta.$$

Since by Property 6.2 $\tilde{U}_n(\varepsilon) \to \tilde{U}(\varepsilon) \in \mathbb{R}$, it follows that

$$U_n(x) \leq C + \varepsilon x \qquad \text{for all } x > -2\beta$$

for some $C > 0$. Since $\tilde{U}_n$ is convex and $U_n$ is nondecreasing, we deduce that

$$U_n(-D\tilde{U}_n(y)) \leq U_n(-D\tilde{U}_n(y_0)) \leq C - \varepsilon D\tilde{U}_n(y_0) \qquad \text{for all } y \geq y_0.$$

Now observe that $D\tilde{U}_n(y_0)$ is bounded uniformly in $n$ by Properties 6.2 and 6.3 together with the closedness of $\{(x,y) : x \in \partial \tilde{U}(y)\}$ (see, e.g., [16]). It follows that there exists some $C > 0$ such that

$$(6.22) \qquad U_n(-D\tilde{U}_n(y)) \leq C \qquad \text{for all } y \geq y_0 \text{ and } n \geq 1.$$

We can now conclude the proof of (6.21). Since $X_n - B = -D\tilde{U}_n(Y_n)$, it follows from Property 6.4 and (6.22) that, on $\{Y_n > 0\}$,

$$I_n \leq C\mathbf{1}_{Y_n \geq y_0} + \{J_n + [(X_n - B) + \beta]Y_n + (B - \beta)Y_n\}\mathbf{1}_{Y_n \leq y_0}$$

$$\leq C + \{J_n + C[1 + \tilde{U}_n(Y_n) - Y_n\beta] + (B - \beta)Y_n\}\mathbf{1}_{Y_n \leq y_0}$$

$$\leq 2C + (C+1)|J_n|,$$

where we used the fact that $B - \beta \leq 0$. It follows that

$$\hat{I}_n^+ \leq \widehat{|J_n|} = \hat{J}_n + 2(\widehat{J_n^-}),$$

where $\widehat{|J_n|}$ (resp. $\widehat{J_n^-}$) is the convex combination in $\text{conv}\{|J_k|, k \geq n\}$ (resp. $\text{conv}\{J_k^-, k \geq n\}$) corresponding to $\hat{I}_n$. Since $\tilde{U}_n(0) = U_n(\infty)$ and $\tilde{U}(Y_n) < \infty$, it follows that $Y_n > 0$, $P$-a.s. whenever $U_n(\infty) = \infty$. Therefore, $I_n^+$ is bounded on $\{Y_n = 0\}$. In view of this, we obtain immediately (6.21) from the uniform integrability of the sequences $(\hat{J}_n)_n$ and $(J_n^-)_n$; see Corollary 6.2 and (6.15). □

COROLLARY 6.3. *Suppose that $L := \inf\{l \geq 0 : U(l) = U(+\infty)\} = +\infty$. Then $Y_* > 0$, $P$-a.s.*

PROOF. The case $U(\infty) = \infty$ is easily treated because it implies that $\tilde{U}(0) = +\infty$. We then concentrate on the case where $U$ is bounded. By Step 1 of the proof of Proposition 3.1 (see the end of Section 5), it follows from the condition $L = +\infty$ that $\partial \tilde{U}(0) = \{-\infty\}$.

Let $P_0 := Y_0 \cdot P$ be an arbitrary measure in $\mathcal{M}^e(S)$. From Lemma 6.4, we have $E[Y_0 X_*] \leq x$. Since $X_* \geq B - 2\beta$, this proves that $E[Y_0(X_*)^+] < \infty$. However, $X_* = +\infty$ on the event set $\{Y^* = 0\}$. Hence $P_0[Y_* = 0] = 0$ and the proof is complete. □

We are now able to complete the proof of Theorem 3.2(ii) and (iii).

LEMMA 6.6. *There exists a random variable $\bar{X}_* \in \mathcal{C}(x)$ that satisfies*

$$xy_* = E\bar{X}_*Y_*, \qquad \bar{X}_* \in B - \partial\tilde{U}(Y_*), \qquad P\text{-a.s.} \quad \text{and} \quad EU(\bar{X}_* - B) = V(x).$$

*Moreover, if $\bar{X}_* \geq 0$, then $\bar{X}_* \in \mathcal{X}_+(x)$.*

PROOF.

STEP 1.   Combining Lemmas 6.4 and 6.5, we see that $\bar{X}_* := X_* \in \mathcal{C}(x)$ and satisfies the announced requirements.

STEP 2.   We now assume that $X_* \geq 0$, $P$-a.s. By Remark 3.5, there exists some $\tilde{X}_* \in \mathcal{X}_+(x)$ such that $\tilde{X}_* \geq X_*$, $P$-a.s. Since $\mathcal{X}_+(x) \subset \mathcal{C}(x)$, we have $E\tilde{X}_*Y_* \leq xy_* = EX_*Y_*$ and therefore $\tilde{X}_* = X_*$ on $\{Y_* > 0\}$. We next consider two cases.

2.1. Assume first that $L := \inf\{l \geq 0 : U(l) = U(+\infty)\} = +\infty$. Then, from Corollary 6.3, $Y_* > 0$, $P$-a.s. It follows that $\tilde{X}_* = X_*$, $P$-a.s. and the requirement of the lemma holds for $\bar{X}_* := X_*$.

2.2. If $L < \infty$, then $Y_*$ may be zero with positive probability. However, since $\tilde{X}_* = X_*$ on $\{Y_* > 0\}$ and $X_* - B \in -\partial\tilde{U}(Y_*)$, we have

$$E[\tilde{X}_*Y_*] = xy_* \quad \text{and} \quad (\tilde{X}_* - B) \wedge L = (X_* - B) \wedge L.$$

Since $U(x) = U(L)$ for $x \geq L$, this proves that

$$\tilde{X}_* \in B - \partial\tilde{U}(Y_*) \quad \text{and} \quad U(\tilde{X}_* - B) = U(X_* - B), \qquad P\text{-a.s.}$$

Hence, the required result holds for $\bar{X}_* := \tilde{X}_*$.

This completes the proof of Theorem 3.2(ii) and (iii).   $\square$

6.4. *Proof of Lemma 6.1.*   The last statement of the lemma follows from a direct application of the la Vallée–Poussin theorem. Let $n$ be a fixed integer in $[1, \infty]$ and consider the two following cases.

CASE 1.   Suppose that $\tilde{U}_n(+\infty) = -\infty$. Then

$$\tilde{U}_n : (0, \infty) \to (-\infty, \tilde{U}_n(0)) \qquad \text{is convex and decreasing.}$$

Observe that this is valid even for the case $n = \infty$, where $\tilde{U}_\infty = \tilde{U}$ is not strictly convex. Let

$$\phi_n : (-\tilde{U}_n(0), +\infty) \to (0, \infty)$$

be the inverse of $-\tilde{U}_n$. By direct computation, we see that for all $(y, Y) \in \mathbf{Y}_+$ with $y > 0$,

$$E[\phi_n(\tilde{U}_n(Y)^-)] = E[\phi_n(\max\{0, -\tilde{U}_n(Y)\})]$$

$$\leq E[\max\{\phi_n(0), Y\}]$$

$$\leq \phi_n(0) + E[Y] \leq \phi_n(0) + y.$$

Recall that $\tilde{U}(0) = U(+\infty) > 0$ by (6.1), so that $\phi_\infty(0) < \infty$. By increase of $(\tilde{U}_n)_n$, we deduce that $(\phi_n)_n$ is increasing and therefore $\phi_n(0) \leq \phi_\infty(0) < \infty$.

It remains to prove that $\lim_{x\to\infty}[\phi_n(x)/x] = \infty$ or, equivalently, by a trivial change of variable,

$$\lim_{y\to+\infty} \frac{y}{-\tilde{U}_n(y)} = \infty. \tag{6.23}$$

Let us consider separately the cases $n = \infty$ and $n < \infty$.

1. If $n = \infty$, then by an easy extension of l'Hôpital's rule to the nonsmooth case, we see that

$$\lim_{y\to\infty} \frac{y}{-\tilde{U}_n(y)} \geq \liminf_{y\to\infty} \inf_{q\in-\partial\tilde{U}(y)} \frac{1}{q} = \liminf_{y\to\infty} \left[\sup_{q\in-\partial\tilde{U}(y)} q\right]^{-1}.$$

Now, recall that $\tilde{U}(\infty) = U(0) = -\infty$, and therefore $\lim_{x\to 0} \inf \partial U(x) = \infty$ and $\lim_{y\to\infty} \sup -\partial\tilde{U}(y) = 0$ by the classical connection between the generalized gradients of $U$ and $\tilde{U}$. This provides (6.23).

2. If $n < \infty$, then by l'Hôpital's rule together with Property 6.3 (with $\beta = 0$), we see that

$$\lim_{y\to+\infty} \frac{y}{-\tilde{U}_n(y)} = \lim_{y\to+\infty} \frac{1}{-n(y - z_n(y))}, \tag{6.24}$$

where $z_n(y)$ is defined in Property 6.1. Now, from the definition of $\tilde{U}$ and $\tilde{U}_n$ together with (6.2), we have

$$U(x) - xz_n(y) + \frac{n}{2}|z_n(y) - y|^2 \leq \tilde{U}(z_n(y)) + \frac{n}{2}|z_n(y) - y|^2$$
$$= \tilde{U}_n(y) \leq \tilde{U}(y)$$

for all $x > 0$. Then

$$\frac{n}{2}|z_n(y) - y|^2 \leq \tilde{U}(y) - U(x) + x(z_n(y) - y) + xy$$

$$\leq \tilde{U}(y) - U(x) + xy + \frac{|x|^2}{n} + \frac{n}{4}|z_n(y) - y|^2,$$

where we used the trivial inequality $ab \leq na^2/4 + b^2/n$. This provides

$$\frac{n}{4}|z_n(y) - y|^2 \leq \tilde{U}(y) - U(x) + xy + \frac{|x|^2}{n}.$$

In particular, taking $x = \hat{x}_y \in -\partial\tilde{U}(y)$, we have $\tilde{U}(y) - U(\hat{x}_y) + y\hat{x}_y = 0$ and

$$\frac{n}{4}|z_n(y) - y|^2 \leq \frac{|\hat{x}_y|^2}{n} \leq \frac{1}{n} \sup_{q\in-\partial\tilde{U}(y)} |q|^2.$$

Since $U(0) = -\infty$, it follows that $\inf\{|p| : p \in \partial U(x)\} \to +\infty$ as $x \searrow 0$ and therefore $\sup\{|q| : q \in -\partial \tilde{U}(y)\} \to 0$ as $y \nearrow \infty$ by the classical connection between the generalized gradients of $U$ and $\tilde{U}$. Hence, the last inequality proves that $n|z_n(y) - y| \to 0$ as $y \nearrow \infty$, and (6.23) follows from (6.24).

CASE 2.  We now consider the case where $\tilde{U}_n(+\infty) > -\infty$. We reduce the problem to that of Case 1 by defining the function

$$\phi_n(z) := \begin{cases} (-\tilde{U}_n)^{-1}(z), & \text{for } -\tilde{U}_n(0) \le z \le -\tilde{U}_n(+\infty), \\ \psi_n(z), & \text{for } z \ge -\tilde{U}_n(+\infty), \end{cases}$$

where $\psi_n$ is chosen so that $\phi_n(x)/x \to +\infty$ as $x \nearrow \infty$. It is immediately checked that the inequality $E[\phi_n(\tilde{U}_n(Y)^-)] \le \phi_n(0) + y$ holds with this definition of $\phi_n$. Finally, arguing as in Case 1, we can choose $(\psi_n)_n$ such that $\phi_n(0) \le \phi_\infty(0) < \infty$.

**7. The asymptotic elasticity conditions.**  In this section we prove Lemma 2.3 which has been used extensively for the proof of our main result.

PROOF OF LEMMA 2.3.

STEP 1.  From the nonincrease of $f$ near zero, we have

$$\mathrm{AE}_0(f) = \limsup_{y \downarrow 0} \sup_{q \in \partial f(y)} \frac{-qy}{f(y)}.$$

This is in agreement with the definition of [5], where Lemma 4.1 states that the asymptotic elasticity condition $\mathrm{AE}_0(f) < \infty$ is equivalent to the existence of $y_0 > 0$ and $\beta > 0$ such that

$$f(\mu y) \le \mu^{-\beta} f(y) \qquad \text{for all } \mu \le 1 \text{ and } y \le y_0.$$

STEP 2.  By a similar argument to Lemma 4.1 in [5], we also obtain a characterization of the asymptotic elasticity condition $\mathrm{AE}_\infty(f) < \infty$ by the existence of $y_1 > 0$ and $\beta > 0$ such that

$$f(\mu y) \le \mu^{\beta} f(y) \qquad \text{for all } \mu \ge 1 \text{ and } y_1 \le y.$$

STEP 3.  Since $f$ is nondecreasing near $+\infty$ and nonincreasing near 0, it follows from Steps 1 and 2 that statement (i) of Lemma 2.3 holds for all $y \in (0, y_0] \cup [y_1, \infty)$ (after possibly changing $y_0$ and $y_1$). Since $f(y) > 0$, the inequality of (i) holds on the interval $(y_0, y_1)$ by a simple compactness argument.

STEP 4. We finally prove (ii). Given $y > 0$, let $q$ be an arbitrary element of $\partial f(y)$. By convexity of $f$ together with (i), we have

$$(\mu - 1)yq \leq f(\mu y) - f(y) \leq (C - 1)f(y)$$

for all $\mu \in [2^{-1}, 2]$. The required result is obtained by taking the values $\mu = 2$ and $\mu = 2^{-1}$. □

B. BOUCHARD
LFA-CREST AND UNIVERSITÉ PARIS 6
15 BOULEVARD GABRIEL PERI
92245 MALAKOFF CEDEX
FRANCE
E-MAIL: bouchard@ccr.jussieu.fr

N. TOUZI
LFA-CREST
15 BOULEVARD GABRIEL PERI
92245 MALAKOFF CEDEX
FRANCE
E-MAIL: touzi@ensae.fr

A. ZEGHAL
CEREMADE
UNIVERSITÉ PARIS DAUPHINE
PLACE DU MARECHAL DE
    LATTRE DE TASIGNY
75016 PARIS
FRANCE
E-MAIL: zeghal@ceremade.dauphine.fr